\documentclass[preprint,12pt]{elsarticle}
\usepackage{amsmath}

\usepackage{amssymb}
\newcommand{\citeaynum}[1]{\citeauthor{#1} (\citeyear{#1}) \cite{#1}}
\usepackage{graphicx}
\usepackage{pdfpages}
\usepackage{float}
\usepackage{caption} 
\usepackage{hyperref}

\usepackage[utf8]{inputenc}
\usepackage[T1]{fontenc}

\journal{}

\begin{document}

\begin{frontmatter}



\title{A new high-order finite-volume advection scheme on spherical Voronoi grids and a comparative study in a mimetic finite-volume moist shallow-water model}


\author[inst1]{Luan F. Santos\corref{cor1}}
\ead{luan.santos@usp.br}

\author[inst2]{Jeferson B. Granjeiro}
\ead{matebram@gmail.com}

\author[inst2]{Pedro S. Peixoto}
\ead{ppeixoto@usp.br}

\cortext[cor1]{Corresponding author.}

\affiliation[inst1]{organization={Institute of Mathematical and Computer Sciences, University of São Paulo},
            addressline={Avenida Trabalhador São-carlense, 400},
            city={São Carlos},
            postcode={13566-590},
            state={SP},
            country={Brazil}}

\affiliation[inst2]{organization={Institute of Mathematics, Statistics and Computer Science, University of São Paulo},
            addressline={Rua do Matão, 1010, Butantã},
            city={São Paulo},
            postcode={05508-090},
            state={SP},
            country={Brazil}}
\begin{abstract}
Spherical centroidal Voronoi tessellations (SCVTs), currently used in numerical weather forecasting models such as the Model for Prediction Across Scales (MPAS), are a type of spherical grid that is highly flexible, allowing the construction of locally refined regions with higher resolution without requiring modifications to the numerical discretization or its implementation. However, the irregularity of SCVT grids makes the construction of robust high-order schemes challenging. In particular, in atmospheric modeling, high-order advection schemes are desirable since they reduce numerical diffusion and improve the representation of fine-scale tracer structures. Therefore, in this work, we propose a new class of high-order advection schemes on the sphere based on the $k$-exact reconstruction approach, extending their successful use on planar domains to the spherical surface. We assess the performance of the proposed method and compare it with existing advection schemes for SCVT grids used in MPAS.  The evaluation includes classical advection test cases on the sphere as well as simulations with a mimetic finite-volume moist shallow-water model, in which the advection scheme is applied to the transport of moisture tracers. Grid-related robustness was investigated using locally refined spherical grids with a local focus on the Andes topography. Our results show that the proposed schemes achieve high-order accuracy in the advection tests, exhibit little sensitivity to grid distortion, and produce comparable results to existing schemes in the moist shallow-water model. Overall, grid robustness is therefore limited to the sensitivity of the discretization of the shallow-water model, irrespective of the advection scheme.
\end{abstract}



\end{frontmatter}


\section{Introduction}
The increasing demand for higher spatial resolutions in global atmospheric dynamical cores has driven the development of scalable, massively parallel numerical schemes and motivated the development of alternative spherical grids for atmospheric modeling \citep{stan2012grids}. These grids differ from traditional latitude–longitude grids, which have several well-known limitations for explicit time-stepping schemes, particularly those related to grid singularities near the poles and restrictive Courant-Friedrichs-Lewy (CFL) stability constraints \citep{will2007dycores}. 

Among alternative global atmospheric grids, cubed-sphere, icosahedral, and pentagonal–hexagonal grids have received particular attention because of their quasi-uniform spacing and suitability for parallel computing \citep{stan2012grids}.
These grids are generated by projecting Platonic solids onto the sphere. Finite-volume methods are well suited to them because they are flexible with respect to grid geometry, preserve local mass, and rely on local flux computations that favor parallelism. Still, challenges remain, including grid imprinting \citep{mouallem2024duogrid,peixoto2013gridimprinting,weller2012gridimprinting} and convergence issues caused by problematic spherical cells \citep{peixoto2013gridimprinting,peixoto2016consitency}.

To improve uniformity and reduce grid-induced errors, Platonic-solid grids are often optimized \citep{tomita2002springdynamics,miura2005quality,heikes2013optimized}.
For instance, Lloyd’s method \citep{du1999scvt,du2003scvt} can be used to generate a spherical centroidal Voronoi tessellation (SCVT).
In SCVT grids, each Voronoi cell’s center of mass coincides with its generator, which helps reduce grid imprinting \citep{peixoto2013gridimprinting}.
Additionally, Lloyd’s method can incorporate a density function to control cell spacing, and SCVT grids allow highly flexible local refinement \citep{engwirda2017jigsaw}.

SCVT grids have been employed in the Model for Prediction Across Scales (MPAS), both in its atmospheric \citep{skamarock2012mpas} and oceanic \citep{ringler2013multiresolgrid} components, together with the mimetic finite-volume/difference TRiSK discretization \citep{thuburn2009TRiSK,ringler2010trsk}.
More recently, MPAS was adopted as the dynamical core of the Model for Ocean-laNd-Atmosphere predictioN (MONAN) \citep{freitas2023monan}, the new community Earth system model developed by the Brazilian National Institute for Space Research (INPE).

Besides its use in MPAS, TRiSK has also been implemented in other dynamical cores \citep{kevlahan2019wavetrisk,dubos2015dynamico}. However, as a low-order method, it may still suffer from grid imprinting and lack of convergence on SCVT grids \citep{peixoto2013gridimprinting,peixoto2016consitency,santos2021andes}. In contrast, the transport scheme in MPAS, proposed by \citet{skamarockgassmann2011} (hereafter SG2011), adopts a higher-order approach.

Advection schemes on pentagonal/hexagonal grids have been widely developed within the finite-volume framework for global atmospheric transport. A key contribution is the second-order upwind method of \citet{miura2007upwind}, based on the swept-area concept for stable scalar transport. It was later extended to second- and fourth-order reconstructions with conservative corrections that preserve cell means \citep{skamarock2010transport,miura2013upwind}. These reconstructions are built on the local tangent plane by least squares, with fluxes evaluated by Gaussian quadrature. High-order extensions to SCVT meshes improved practical accuracy despite reduced formal convergence from mesh irregularity \citep{skamarockgassmann2011}. Later fifth- and sixth-order variants reduced absolute errors at higher cost \citep{zhang2018extending}, while Gaussian-quadrature upwind-biased schemes improved stability and robustness \citep{subich2018higher}.

Despite these developments, constructing accurate and efficient high-order finite-volume schemes on spherical unstructured meshes remains challenging. In this context, we investigate the extension to the sphere of the high-order reconstruction approach developed by Ollivier-Gooch and collaborators, successfully applied to two-dimensional unstructured meshes \citep{ollivier1997high,olliviergooch2002advdifeq,ollivier2009,michalak2009euler,nejat2008gmres,jalali2014difflux,hoshyari2018}. This methodology builds on the work of \citet{barth1990quadratic} and \citet{barth1993kexact}, who introduced the $k$-exact reconstruction technique in computational fluid dynamics for aerospace applications. In this approach, the solution in each control volume is approximated by a polynomial that preserves local conservation over neighboring control volumes.

Building on this framework, we develop a high-order advection scheme on the sphere by extending the methodology of Ollivier-Gooch and collaborators to the spherical surface using Voronoi meshes.
The method employs local tangent planes onto which the variables are projected, as in \citet{skamarock2010transport} and \citet{skamarockgassmann2011}.
We then propose second-, third-, and fourth-order schemes by increasing the degree of the polynomial reconstruction and the number of quadrature points.

We first assess the proposed schemes using standard spherical advection test cases \citep{nair2010class}. We then consider the moist shallow-water model of \citet{zerroukat2015moistSWE}, in which the new advection schemes are applied to the transport of temperature and moisture tracers, while the remaining dynamics are discretized as in \citet{santos2021andes} using the TRiSK mimetic finite-volume/difference scheme \citep{thuburn2009TRiSK,ringler2010trsk}.
In all simulations, we compare the proposed schemes with the SG2011 schemes \citep{skamarockgassmann2011} on both quasi-uniform SCVT meshes and topography-based locally refined grids developed by \citet{santos2021andes}.

The results show that the proposed schemes achieve high-order accuracy in standard spherical advection benchmarks and are insensitive to grid distortion. In contrast, the centered SG2011 schemes are more sensitive to grid distortion, whereas the upwind-biased third-order SG2011 variant is more robust. The proposed methods also perform well on locally refined SCVT meshes. Among them, the fourth-order scheme provides the highest accuracy in the advection tests, although at a higher computational cost.

In the moist shallow-water simulations, the proposed schemes yield results comparable to those of SG2011. However, all schemes remain influenced by the underlying TRiSK discretization, even when higher-order advection is used. This indicates that advances in the shallow-water dynamical discretization itself are still needed, in addition to more accurate tracer transport. Such improvements are also relevant for MPAS-based models, such as MONAN.

This paper is organized as follows. Section~\ref{sec:model-description} provides a description of the model, including the SCVT grids used in this work, the finite-volume discretization of the advection equation on the SCVT grid, an overview of the schemes proposed by SG2011, and an introduction to the new high-order advection schemes, as well as the mimetic finite-volume moist shallow-water model.
Section~\ref{sec:num-exp} presents the numerical experiments for both the advection and moist shallow-water equations on the sphere.
Finally, conclusions are drawn in Section~\ref{sec:conclusions}.

\section{Model description}
\label{sec:model-description}
In this section, we introduce the discrete framework used in this work. We first present the quasi-uniform and locally refined Voronoi grids, together with the notation used throughout. We then describe the finite-volume discretization of the advection equation on the sphere, review the SG2011 method, and introduce the proposed scheme. Finally, we present the moist shallow-water model and its mimetic finite-volume discretization, in which the advection schemes are applied to the model tracers.

\subsection{Grids}
\label{sec:grids}
In this work, we consider Spherical Centroidal Voronoi Tessellation (SCVT) grids. 
We briefly explain them here; for a complete survey on SCVT grids, we refer to \citet{ju2011scvt}.

\begin{figure}[!htb]
    \centering
    \includegraphics[width=\textwidth]{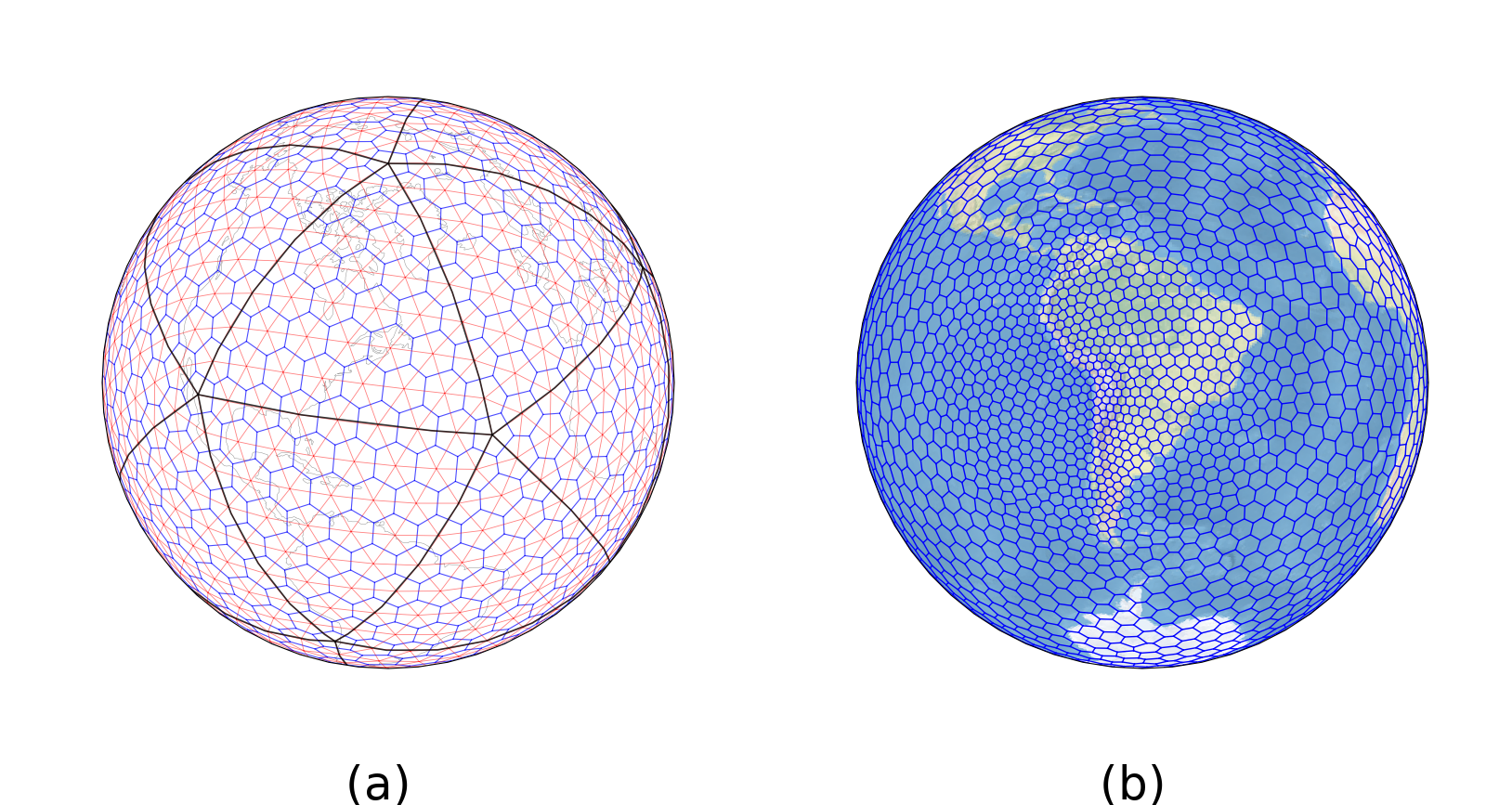}
    \caption{Illustration of SCVT grids used in this work. (a) Spherical grid composed of pentagonal and hexagonal cells (grid level $=3$). Black lines show the base icosahedron, red lines the refined Delaunay triangulation, and blue lines the optimized SCVT Voronoi cells. (b) Locally refined SCVT grid over the Andes mountains (grid level $=4$).}
    \label{fig:scvt_grid}
\end{figure}

To introduce SCVT grids, we start with a pentagonal/hexagonal grid, which is related to the icosahedral grid in the graph-theoretic sense, where they are dual to each other. 
The icosahedral grid is generated from an initial base icosahedron, as illustrated in Figure~\ref{fig:scvt_grid}a. 
This grid is then refined, producing the Delaunay triangulation shown by the red lines in Figure~\ref{fig:scvt_grid}a. 
The pentagonal/hexagonal grid is obtained by connecting the circumcenters of these triangles, 
which are referred to as Voronoi generators. The resulting grid is known as a Voronoi tessellation. The level of refinement of the icosahedron is referred to as the grid level. For a grid level $l$, the total number of Voronoi cells is 
\begin{equation*}
    N_c = 10 \times 2^l + 2.
\end{equation*}

Voronoi tessellations can be optimized with Lloyd’s method \citep{du1999scvt,du2003scvt} so that the center of mass of each cell, with respect to a given density function, coincides with its generator, yielding the spherical centroidal Voronoi tessellations (SCVTs) used in this work. With a suitable density function, Lloyd’s method can also generate locally refined grids. We also consider such grids, in particular the topography-based locally refined grids over the Andes developed by \citet{santos2021andes}, illustrated in Figure~\ref{fig:scvt_grid}b. These grids are generated from a density function based on the Andes topography, producing a resolution about three times finer over the Andes than on the coarse grid. The resolution then transitions smoothly to the surrounding circular grid over South America and, from there, to the coarse grid (Figure~\ref{fig:scvt_grid}b).

Now that we have introduced SCVT grids, we introduce some notation that will be useful for the remainder of this work. 
Inspired by \citet{ringler2010trsk}, we adopt the notation for Voronoi meshes summarized in Table~\ref{tab:notation}.

\begin{table}[!ht]
\centering
\begin{tabular}{ll}
\hline
Symbol & Meaning \\ 
\hline
$i$ & Index of a Voronoi cell center \\ 
$N_c$ & Total number of Voronoi cells \\ 
$\boldsymbol{x}_i$ & Coordinates of cell center $i$ \\ 
$\Omega_i$ & Voronoi cell (control volume) associated with center $i$ \\ 
$|\Omega_i|$ & Area of cell $\Omega_i$ \\ 
$NB(i)$ & Set of Voronoi cells neighboring cell $i$ \\ 
$e$ & Index of a Voronoi edge \\ 
$N_e$ & Total number of edges \\ 
$\boldsymbol{x}_e$ & Coordinates of edge $e$  point\\ 
$\boldsymbol{x}_e^{\text{mid}}$ & Midpoint coordinates of edge $e$\\
$\boldsymbol{n}_e$ & Unit normal vector at edge $e$ \\ 
$n_{e,i}$ & Orientation factor:
$n_{e,i}=
\begin{cases}
+1, & \text{if } \boldsymbol{n}_e \text{ points outward from } \Omega_i,\\
-1, & \text{if } \boldsymbol{n}_e \text{ points inward toward } \Omega_i
\end{cases}$

\\
$EC(i)$ & Set of edges of cell $i$ \\ 
$CE(e)$ & Set of cells adjacent to edge $e$ \\ 
$\partial \Omega_i$ & Boundary of cell $\Omega_i$ \\ 
$\Gamma_e$ & Individual edge of a cell \\ 
$|\Gamma_e|$ & Length of edge $e$ \\ 
\hline
\end{tabular}
\caption{Summary of notation for Voronoi cells, edges, neighbors, and control volumes.}
\label{tab:notation}
\end{table}

Whenever we write $i,j \in CE(e)$, we assume $i \neq j$, where $CE(e)$ denotes the set of cells adjacent to edge $e$. 

We emphasize that $\boldsymbol{x}_e^{\text{mid}}$ (the Voronoi edge midpoint) and $\boldsymbol{x}_e$ (the Delaunay triangle midpoint) are distinct on SCVT grids, as noted by \citet{peixoto2016consitency}. This distinction affects both the linear interpolation of scalar fields at edge points from surrounding cell averages and the midpoint-rule integration used for edge integrals.

\subsection{Finite-volume discretization of the advection equation on the sphere}
\label{sec:fv-adv}
Now that we have described the grids considered in this work, we introduce the finite-volume discretization of the advection equation. We first present the general discretization on the sphere, where the scheme reduces to flux computations based on reconstructing the scalar field from cell averages. We then review the reconstruction of \citet{skamarockgassmann2011} and introduce our new scheme, which extends the planar approach of \citet{olliviergooch2002advdifeq} to the sphere. The notation from Table~\ref{tab:notation} is used throughout.

\subsubsection{General finite-volume discretization}
\label{sec:fvgeneral}
The advection equation on the sphere is given by
\begin{equation}
    \frac{\partial \phi}{\partial t} = - \nabla \cdot (\phi \, \boldsymbol{u}),
    \label{eq:adv}
\end{equation}
where $\phi=\phi(\boldsymbol{x},t)$ is a scalar tracer density, and $\boldsymbol{u}(\boldsymbol{x},t)$ is the fluid velocity vector at a point $\boldsymbol{x}$ on the sphere at time $t$. The cell-averaged tracer is defined as
\begin{equation}
    \bar{\phi}_i(t) = \frac{1}{|\Omega_i|} \int_{\Omega_i} \phi(\boldsymbol{x},t) \, d\Omega_i.
\end{equation}

Integrating Eq.~\eqref{eq:adv} over the control volume $\Omega_i$ (a Voronoi cell) and applying the divergence theorem gives the integro-differential form
\begin{equation}
    \frac{d \bar{\phi}_i}{dt} = - \frac{1}{|\Omega_i|} \int_{\partial \Omega_i} \phi \, \boldsymbol{u} \cdot \boldsymbol{n} \, dl
    = - \frac{1}{|\Omega_i|} \sum_{e \in EC(i)} n_{e,i} \int_{\Gamma_e} \phi \, \boldsymbol{u} \cdot \boldsymbol{n}_e \, dl.
    \label{eq:fv-adv}
\end{equation}
Here, $\bar{\phi}_i$ denotes the cell-averaged tracer concentration, and $\boldsymbol{n}$ is the outward unit normal vector to $\partial \Omega_i$.

A finite-volume discretization of the advection equation then reduces to finding a numerical flux ${F}_e$ that approximates the flux across an edge $\Gamma_e$, i.e.,
\begin{equation}
    {F}_e \approx \int_{\Gamma_e} \phi \, \boldsymbol{u} \cdot \boldsymbol{n}_e \, dl.
\end{equation}
using the cell-averaged values $\bar{\phi}_i$ from the neighboring cells.
Once the numerical flux is known, the scheme reads
\begin{equation}
    \frac{d \bar{\phi}_i}{dt} = - \frac{1}{|\Omega_i|} \sum_{e \in EC(i)} n_{e,i} F_e ,
    \label{eq:fv-adv-discrete}
\end{equation}
which is a system of ordinary differential equations (ODEs) for the cell-averaged values $\bar{\phi}_i$ on each control volume.
This ODE system is then integrated in time using the three-stage Runge-Kutta scheme described in \citet{skamarockgassmann2011} and detailed in Appendix A.

In addition, following \citet{skamarockgassmann2011}, we adopt the flux-corrected transport (FCT) scheme of \citet{zalesak1979fct} to prevent nonphysical negative values and overshoots or undershoots. It modifies only the last Runge–Kutta stage by blending a first-order upwind flux with the higher-order flux $F_e$ and limiting the flux when a new extremum is detected. Further details are given in Appendix B.

\subsubsection{SG2011 scheme}
\label{sec:fv-sg}
For simplicity, we refer to the schemes of SG2011 as SG schemes and omit explicit time dependence, since only the spatial discretization is considered here. The SG schemes compute the flux across an edge using the midpoint rule for integration.
\begin{equation}
\label{eq:midpoint-flux}
 \int_{\Gamma_e} \phi \, \boldsymbol{u} \cdot \boldsymbol{n}_e \, dl 
 = \phi(\boldsymbol{x}_e^{\text{mid}}) \, u_e \, |\Gamma_e| + \mathcal{O}(|\Gamma_e|^2),
\end{equation}
where \(u_e = \boldsymbol{u}(\boldsymbol{x}_e^{\text{mid}}) \cdot \boldsymbol{n}_e\) is the normal velocity at the midpoint.

The tracer value at the edge midpoint, $\phi(\boldsymbol{x}_e^{\text{mid}})$, must be approximated.  
We denote this approximation by $\phi_e^{\text{mid}}$, which then defines the numerical flux as
\begin{equation}
\label{eq:sg-numerical-flux}
 F_e = \phi_e^{\text{mid}} \, u_e \, |\Gamma_e|.
\end{equation}

SG2011 proposed a second-order scheme, denoted {SG2}, in which the midpoint tracer is defined as the arithmetic average of the cell-averaged values of the two cells sharing edge \(e\):\begin{equation}
    \phi_e^{\text{mid}} = \frac{\bar{\phi}_i + \bar{\phi}_j}{2}, \quad i,j \in CE(e),
\end{equation}
where \(CE(e)\) denotes the set of cells adjacent to edge \(e\).
In SG2, the displacement between the edge midpoint $\boldsymbol{x}_e^{\text{mid}}$  and the actual edge point  $\boldsymbol{x}_e$ reduces the order of linear interpolation from two to one \citep{peixoto2016consitency}.

To achieve higher-order accuracy, SG2011 extended one-dimensional finite-volume reconstructions to the sphere. To illustrate the idea, consider the one-dimensional case, where edge \(e\) is denoted by \(i+\tfrac{1}{2}\). The control-volume interface value is reconstructed as
\begin{equation}
\label{eq:ppm-recon}
    \phi_{i+\frac{1}{2}} = \frac{7}{12}(\bar{\phi}_i +\bar{\phi}_{i+1}) - \frac{1}{12}(\bar{\phi}_{i-1} +\bar{\phi}_{i+2}).
\end{equation}
This formula is fourth-order accurate on uniform one-dimensional grids and is used in several schemes based on the piecewise parabolic method (PPM) \citep{colella1984ppm} to reconstruct the advected field  
(e.g., \citep{putman2007cubedsphere, santos2025div, bendall2025swift}).

To obtain a more compact form suitable for generalization, SG2011 rewrote Equation ~\eqref{eq:ppm-recon} as
\begin{equation}
\label{eq:ppm-recon-alt}
    \phi_{i+\frac{1}{2}} = \frac{1}{2}(\bar{\phi}_i + \bar{\phi}_{i+1}) - \frac{1}{12} \left( \delta_x^2 \bar{\phi}_i + \delta_x^2 \bar{\phi}_{i+1} \right),
\end{equation}
where the second-difference operator is defined by
\begin{equation}
    \delta_x^2 \bar{\phi}_i = \bar{\phi}_{i+1} - 2\bar{\phi}_i + \bar{\phi}_{i-1}.
\end{equation}

To extend this one-dimensional formulation to the sphere, the derivative along the Cartesian coordinate $x$ is replaced by the directional derivative along the unit normal $\boldsymbol{n}_e$ to the edge. SG2011 extends Equation~\eqref{eq:ppm-recon-alt} to the sphere as
\begin{equation}
\label{eq:recon-sg4}
    \phi_e^{\text{mid}} = \frac{\bar{\phi}_i + \bar{\phi}_j}{2}
    - \frac{\Delta x_e^2}{12} \left( D^2_{\boldsymbol{n_e}} \phi_i + D^2_{\boldsymbol{n_e}} \phi_j \right),
    \quad i,j \in CE(e),
\end{equation}
where $\Delta x_e$ is the geodesic distance between the control volume centers $\boldsymbol{x}_i$ and $\boldsymbol{x}_j$, and $D^2_{\boldsymbol{n_e}}\phi_i$ and $D^2_{\boldsymbol{n_e}}\phi_j$ are the second-order directional derivatives of $\phi$ at $\boldsymbol{x}_i$ and $\boldsymbol{x}_j$ along $n_{e,i}\boldsymbol{n}_e$ and $n_{e,j}\boldsymbol{n}_e$, respectively.
Because this is a direct extension of a fourth-order one-dimensional scheme, the resulting scheme is denoted by SG4.

To complete the reconstruction, the second-order derivatives of $\phi$ required in Equation~\eqref{eq:recon-sg4} must be computed. In SG2011, this is achieved by fitting a second-order polynomial to the values of $\phi$ at the neighboring control volume centers. 

For each cell center $i$, a quadratic local polynomial is constructed by radially projecting the surrounding cell centers $k \in \mathcal{S}_i$ onto the tangent plane at $\boldsymbol{x}_i$, with coordinates $(x,y)$ centered at $\boldsymbol{x}_i$, where
\begin{equation}
\label{eq:sg4-stencil}
    \mathcal{S}_i = NB(i) \cup \{i\}.
\end{equation}

\begin{figure}[!htb]
    \centering
        \includegraphics[width=1.\textwidth]{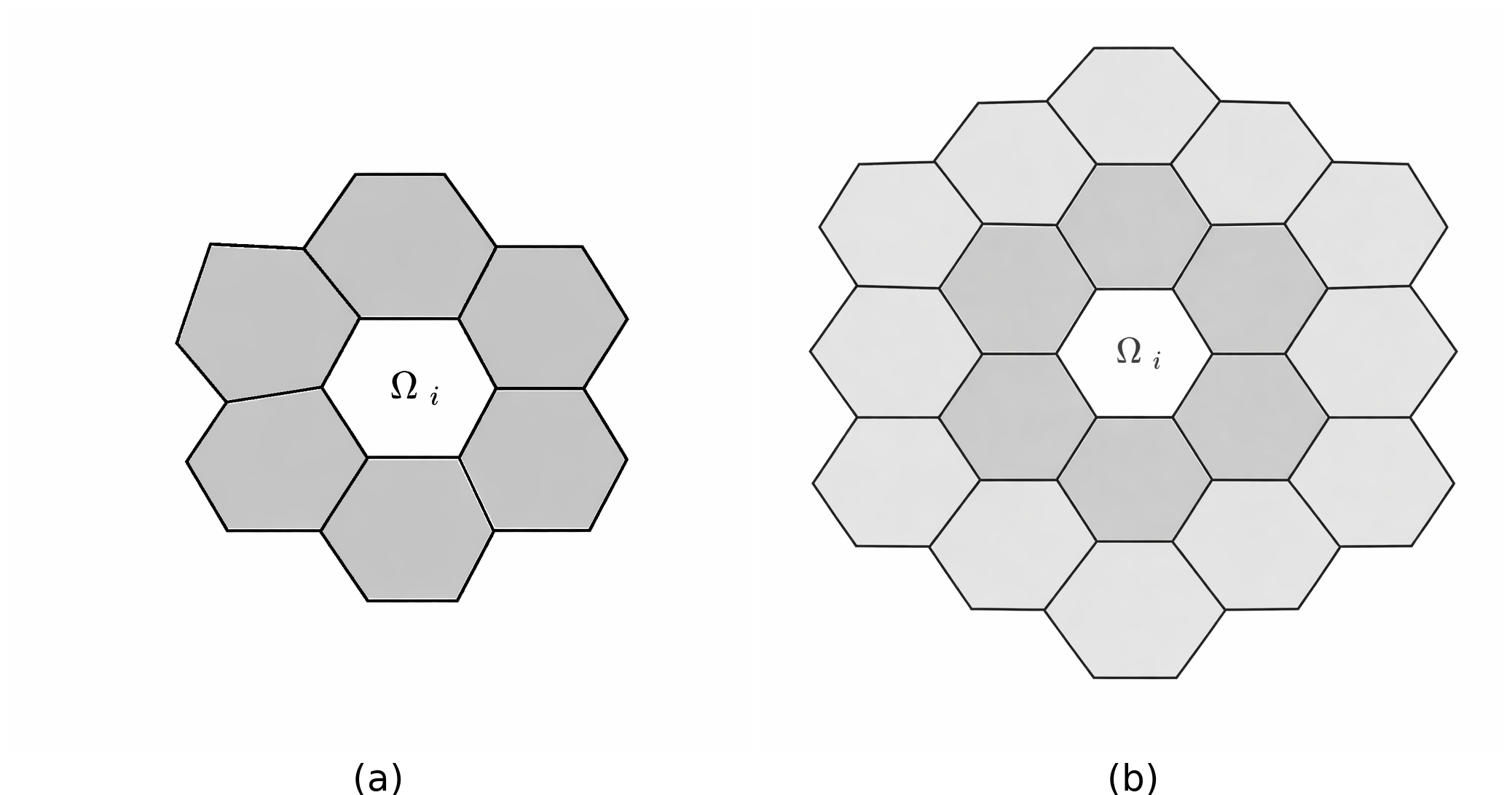}
\caption{Illustration of the stencils using first-level neighbors (a) and first- and second-level neighbors (b), used to construct the reconstruction polynomials.}    \label{fig:scvt_stencil}
\end{figure}

In other words, $\mathcal{S}_i$ contains the first neighbors of cell $i$, as illustrated in Fig.~\ref{fig:scvt_stencil}.

The reconstruction polynomial is then written as
\begin{equation}
\label{eq-sg-poly}
    {\phi}_i^R(x,y) = c_0 + c_1 x + c_2 y + c_3 x^2 + c_4 xy + c_5 y^2,
\end{equation}
with coefficients $\{c_m\}$ determined by enforcing
\begin{equation}
\label{eq-sg-poly-system}
    {\phi}_i^R(x_k, y_k) = \overline{\phi}_k, \quad k \in \mathcal{S}_i,
\end{equation}
where $(x_k, y_k)$ are the projected coordinates of the control volume centers on the tangent plane at $\boldsymbol{x}_i$. 
Thus, $\mathcal{S}_i$ defines the stencil used to fit the polynomial.

The corresponding linear system in matrix form is
\begin{equation}
\label{eq:sg-linear-system-matrix}
\begin{bmatrix}
x_{k_1} & y_{k_1} & x_{k_1}^2 & x_{k_1} y_{k_1} & y_{k_1}^2 \\
x_{k_2} & y_{k_2} & x_{k_2}^2 & x_{k_2} y_{k_2} & y_{k_2}^2 \\
\vdots & \vdots & \vdots & \vdots & \vdots \\
x_{k_N} & y_{k_N} & x_{k_N}^2 & x_{k_N} y_{k_N} & y_{k_N}^2
\end{bmatrix}
\begin{bmatrix}
c_1 \\ c_2 \\ c_3 \\ c_4 \\ c_5
\end{bmatrix}
=
\begin{bmatrix}
\overline{\phi}_{k_1} - \overline{\phi}_i \\
\overline{\phi}_{k_2} - \overline{\phi}_i \\
\vdots \\
\overline{\phi}_{k_N} - \overline{\phi}_i
\end{bmatrix},
\end{equation}
where $k_1, \dots, k_N \in \mathcal{S}_i \setminus \{i\}$, $N = |\mathcal{S}_i|-1$, and $c_0 = \bar{\phi}_i$. 

This generally overdetermined system is solved in the least-squares sense, following SG2011. Its pseudoinverse is precomputed and stored for efficient flux evaluation. In this work, we compute it via Singular Value Decomposition (SVD) using LAPACK.

Thus, once the polynomial coefficients are obtained, the reconstruction polynomial is used to evaluate the second-order derivatives along the projected normal direction and then recover the scalar field at the edge through Equation~\eqref{eq:recon-sg4}.

More concretely, the second-order derivative along the projected normal direction is computed as
\begin{equation}
\label{eq-sg-normal-d2x}
D^2_{\boldsymbol{n_e}}\phi_i = \boldsymbol{\hat{n}}_e^T \, ( \mathrm{Hess}\, \phi_i^R ) \, \boldsymbol{\hat{n}}_e,
\end{equation}
where $\boldsymbol{\hat{n}}_e$ is the projection of the normal vector $\boldsymbol{n}_e$ onto the tangent plane at $\boldsymbol{x}_i$, 
the superscript \(T\) denotes the matrix transpose, and $\mathrm{Hess}\, \phi_i^R$ is the Hessian matrix of the reconstruction polynomial at cell $i$ (Equation~\eqref{eq:recon-sg4}), which is a $2\times 2$ matrix in this case.

Finally, SG2011 proposed a third-order scheme, denoted as {SG3}, by extending the one-dimensional upwind-based reconstruction of \citet{wicker2002timeSplitting} to the sphere. The edge value is computed as
\begin{align}
\label{eq:recon-sg3}
    \phi_e^{\text{mid}} &= \frac{\bar{\phi}_i + \bar{\phi}_j}{2}
    - \frac{\Delta x_e^2}{12} \left( D^2_{\boldsymbol{n_e}} \phi_i + D^2_{\boldsymbol{n_e}} \phi_j \right) \notag \\
    &\quad + \mathrm{sign}(n_{e,i}u_e)\,\beta\,\frac{\Delta x_e^2}{12} 
      \left( D^2_{\boldsymbol{n_e}} \phi_j - D^2_{\boldsymbol{n_e}} \phi_i \right),
    \qquad i,j \in CE(e),
\end{align}
where $\beta$ is a tunable parameter that controls the upwind bias of the reconstruction. The unit normal vector $\boldsymbol{n}_e$ is assumed to point from the center of control volume $i$ toward $j$.
Following SG2011, we adopt $\beta = 0.25$ or $\beta = 1$ in our simulations.

Unlike the centered SG2 and SG4 schemes, SG3 includes an upwind correction through the term proportional to $\mathrm{sign}(n_{e,i}u_e)\,\beta$. Setting $\beta=0$ recovers SG4, while $\beta\neq 0$ introduces an upwind bias that can improve numerical stability.

Finally, Eq.~\eqref{eq:recon-sg3} can be extended to higher orders using higher-order one-dimensional analogs, such as the fifth- and sixth-order formulations of \citet{zhang2018extending}. However, as noted by \citet{zhang2018extending}, these schemes do not attain on the sphere the same formal order of accuracy as in one-dimensional settings.

In fact, SG3 and SG4 in \citet{skamarockgassmann2011} exhibit second-order convergence in the standard cosine-bell solid-body rotation test with zonal wind. Although they do not attain their nominal convergence rates, they yield smaller overall errors than lower-order schemes, a behavior also observed in the higher-order schemes of \citet{zhang2018extending}.

\subsubsection{The proposed scheme}
\label{sec:fv-og}
Our goal in this section is to introduce the scheme we propose, which extends the high-order finite-volume methodology of Ollivier-Gooch et al. \citep{ollivier1997high,olliviergooch2002advdifeq,ollivier2009,michalak2009euler,nejat2008gmres,jalali2014difflux,hoshyari2018} to the sphere. For simplicity, we shall refer to these methods as the OG schemes.

The proposed scheme aims to achieve a nominal order that matches the actual convergence rate, which, as pointed out in the previous section, is not the case for SG-type schemes. This is accomplished by employing Gaussian quadrature to evaluate the flux integral in Equation~\eqref{eq:sg-numerical-flux}, together with higher-order polynomial reconstruction.

Furthermore, we adopt the OG least-squares procedure, which guarantees that the reconstruction polynomial conserves the mass of the target control volume and preserves local mass consistency across stencil neighbors. This contrasts with the SG schemes, where the reconstruction polynomial is constrained to match pointwise the cell-averaged values in the stencil. Based on \citet{barth1990quadratic}, the OG reconstruction ensures local mass conservation on unstructured grids and is $k$-exact, that is, exact for polynomials of degree $k$. Preserving local mass is desirable in Godunov-type schemes \citep{godunov1959,vanleer1977,colella1984ppm}, as it is expected to improve accuracy and physical consistency \citep{lauritzen2011}. Details of the proposed scheme are given below.

In SG2, the numerical flux is evaluated using the midpoint rule, as given in Equation~\eqref{eq:sg-numerical-flux}. 
As proposed by \citet{barth1990quadratic}, a more accurate approach is to approximate the flux integral over an edge using Gaussian quadrature:
\begin{equation}
\label{eq:og-numerical-flux}
    F_e = \sum_{l=1}^{m} w_l\, \phi_{e,l}\, u_{e,l},
\end{equation}
where $w_l$ are the quadrature weights, $\phi_{e,l}$ is the reconstructed scalar value at the $l$-th quadrature point along edge $\Gamma_e$, and $u_{e,l}$ is the corresponding normal velocity component. In the moist shallow-water model, discussed later, the edge normal velocity $u_e$ is prognostic, so $u_{e,l}$ must be reconstructed from it.


Next, we reconstruct the scalar field at the Gaussian quadrature points from the cell-averaged values $\bar{\phi}_i$ at the control volume centers. This is done by polynomial fitting on a tangent-plane projection, as in SG2011. On the tangent plane, we apply the reconstruction of \citet{olliviergooch2002advdifeq}, which preserves the local mass of each control volume in the stencil.

On the tangent plane, we consider a reconstruction polynomial of a given degree, which can be expressed as
\begin{align}
\label{eq-og-poly}
{\phi}_i^R(x,y) 
&= c_0 
  + c_1 (x - x_i) 
  + c_2 (y - y_i) \nonumber\\
&\quad
  + c_3 (x - x_i)^2 
  + c_4 (x - x_i)(y - y_i) 
  + c_5 (y - y_i)^2 
  + \cdots,
\end{align}
where the coefficients $\{c_m\}$ are determined by enforcing that the polynomial preserves the cell-averaged values over the stencil $\mathcal{S}_i$
\begin{equation}
\label{eq-og-poly-system}
    \frac{1}{|\Omega_k|}\int_{\Omega_k} {\phi}_i^R(x, y)\, d{\Omega_k} = \overline{\phi}_k, 
    \quad k \in \mathcal{S}_i.
\end{equation}

From Equation~\eqref{eq-og-poly-system}, the corresponding linear system in matrix form can be written as
\begin{equation}
\label{eq:og-linear-system-matrix}
\begin{bmatrix}
1 & \overline{x}_{i} & \overline{y}_{i} & \overline{x_i^2} & \overline{xy_i} & \overline{y_i^2} & \cdots \\[4pt]
w_{k_1} & w_{k_1}\,\widehat{x}_{k_1} & w_{k_1}\,\widehat{y}_{k_1} & 
w_{k_1}\,\widehat{x_{k_1}^2} & 
w_{k_1}\,\widehat{\smash[b]{(x y_{k_1})}} & 
w_{k_1}\,\widehat{y_{k_1}^2} & \cdots \\[4pt]
w_{k_2} & w_{k_2}\,\widehat{x}_{k_2} & w_{k_2}\,\widehat{y}_{k_2} & 
w_{k_2}\,\widehat{x_{k_2}^2} & 
w_{k_2}\,\widehat{\smash[b]{(x y_{k_2})}} & 
w_{k_2}\,\widehat{y_{k_2}^2} & \cdots \\[4pt]
\vdots & \vdots & \vdots & \vdots & \vdots & \vdots & \\[4pt]
w_{k_N} & w_{k_N}\,\widehat{x}_{k_N} & w_{k_N}\,\widehat{y}_{k_N} & 
w_{k_N}\,\widehat{x_{k_N}^2} & 
w_{k_N}\,\widehat{\smash[b]{(x y_{k_N})}} & 
w_{k_N}\,\widehat{y_{k_N}^2} & \cdots ,
\end{bmatrix}
\begin{bmatrix}
c_0 \\ c_1 \\ c_2 \\ c_3 \\ c_4 \\ c_5 \\ \vdots
\end{bmatrix}
=
\begin{bmatrix}
\overline{\phi}_{i} \\
w_{k_1}\,\overline{\phi}_{k_1} \\
w_{k_2}\,\overline{\phi}_{k_2} \\
\vdots \\
w_{k_N}\,\overline{\phi}_{k_N}
\end{bmatrix},
\end{equation}
where $k_1, \dots, k_N \in \mathcal{S}_i \setminus \{i\}$ and $N = |\mathcal{S}_i|-1$, and the weights are defined as
\begin{equation}
    w_{k_j} = \frac{1}{{(x_{k_j}-x_i)^2 + (y_{k_j}-y_i)^2}},
\end{equation}
which assigns larger weights to the cells whose centers are closer to the reference cell $i$.
The terms $\overline{x_i^m y_i^n}$ and $\widehat{x_{k_j}^m y_{k_j}^n}$ denote the cell-averaged monomials over $\Omega_i$ and $\Omega_{k_j}$, respectively, that is,
\begin{align}
\overline{x_i^m y_i^n} &= \frac{1}{|\Omega_i|} 
\int_{\Omega_i} (x - x_i)^m (y - y_i)^n \, d\Omega_i, 
\label{eq:og-monomial-i}\\[4pt]
\widehat{x_{k_j}^m y_{k_j}^n} &= \frac{1}{|\Omega_{k_j}|} 
\int_{\Omega_{k_j}} (x - x_i)^m (y - y_i)^n \, d\Omega_{k_j}.
\label{eq:og-monomial-k}
\end{align}
These monomial integrals are converted into boundary line integrals via the Gauss divergence theorem and evaluated exactly by Gaussian quadrature; see \citet{olliviergooch2002advdifeq}. 
The resulting system is overdetermined and solved in the least-squares sense after eliminating the first column, associated with the constant term $c_0$, by Gaussian elimination as in \citet{olliviergooch2002advdifeq}.

Once the reconstruction polynomial coefficients are computed, the scalar field at the Gaussian quadrature points on edge $\Gamma_e$ is evaluated with the upwind-biased polynomial. Let $i^{\mathrm{UPW}}(e)$ denote the upwind control volume associated with edge $e$, defined by
\begin{equation}
    i^{\mathrm{UPW}}(e) =
    \begin{cases}
        i, & \text{if } n_{e,i} u_e \ge 0, \\
        j, & \text{otherwise},
    \end{cases}
\end{equation}
where $i,j \in CE(e)$.

The scalar field at the quadrature points, needed to compute the numerical flux in Equation~\eqref{eq:og-numerical-flux}, is obtained from the upwind reconstruction polynomial:
\begin{equation}
    \phi_{e,l} = {\phi}_{i^{\mathrm{UPW}}(e)}^R(x_{e,l},y_{e,l}), 
    \quad l = 1, \dots, m,
\end{equation}
where $(x_{e,l},y_{e,l})$ are the coordinates of the $l$-th quadrature point projected onto the tangent plane of the upwind cell.

We consider Gaussian quadrature rules with $m=1$ and $m=2$ points.
For $m=1$, the quadrature reduces to the midpoint rule and employs a second-order (linear) polynomial reconstruction, yielding the OG2 scheme. In contrast to SG2, OG2 preserves second-order accuracy even when using the midpoint rule. Its stencil, $\mathcal{S}_i^{\mathrm{OG2}}$, is defined as in Equation~\eqref{eq:sg4-stencil} and consists of the reference control volume and its immediate neighbors.

For $m=2$, we employ third- (quadratic) and fourth-order (cubic) reconstructions, yielding the OG3 and OG4 schemes. Both share the stencil
\begin{equation}
    \mathcal{S}_i^{\mathrm{OG3}} = \mathcal{S}_i^{\mathrm{OG4}} 
    = \bigcup_{k \in \mathcal{S}_i^{\mathrm{OG2}}} \left( NB(k) \cup \{k\} \right),
\end{equation}
which includes the reference cell, its immediate neighbors, and their neighbors (second-order neighbors; see Fig.~\ref{fig:scvt_stencil}), providing sufficient points for least-squares reconstruction.

Unlike the centered SG2 and SG4 schemes, the OG schemes are upwind-biased. 
Both SG3, SG4 and OG3 employ quadratic polynomial reconstructions; however, OG3 relies on a larger stencil, as suggested by \citet{olliviergooch2002advdifeq}.

\subsubsection{Remarks on the computational cost}
Among the SG and OG schemes, SG2 is the least expensive, as it does not require least-squares reconstruction, followed by OG2. In the pure advection experiments of Section~\ref{sec:num-exp}, OG3 was about 10\% more expensive than SG3 on grid levels 6 and 7, while OG4 added only about 1\% relative to OG3. Since OG4 uses the same stencil size as OG3, its additional cost is minimal, making it highly competitive with SG3 given its robustness and accuracy gains, as shown in Section~\ref{sec:num-exp}.

We point out that our codes are written in Fortran with OpenMP shared-memory parallelism, but the experiments were run on research laptops rather than HPC systems, so the results should be viewed as rough estimates. Since the OG methods are more floating-point intensive than SG while requiring similar stencil-based memory access, their relative cost could likely be reduced through further code optimizations, such as improved preprocessing and caching, which were not explored here. A more detailed performance analysis is left for future work, preferably in an MPI setting within the MPAS and MONAN frameworks.

\subsection{Finite-volume discretization of the moist shallow-water equations on the sphere}
\label{sec:mswm}
Besides the linear advection equation, we also assess the SG and OG schemes with the moist shallow-water model of \citet{zerroukat2015moistSWE}, introduced here. This model extends the nonlinear shallow-water equations with simplified moist processes, including condensation and rain formation, through a three-state physics scheme. As a two-dimensional framework, it is much less expensive than full three-dimensional models and is therefore well suited for testing numerical schemes in the context of the physics–dynamics coupling problem \citep{gross:2018}.

The moist shallow-water model introduced by \citet{zerroukat2015moistSWE} has been investigated in some studies, including \citet{ferguson2019amr}, which aimed to assess adaptive mesh refinement; \citet{santos2021andes}, which focused on the impact of topography-based local refinement in SCVT grids using a finite-volume discretization; and more recently by \citet{hartney2025moistSWE}, who employed a compatible finite element discretization.
In the latter, the model of \citet{zerroukat2015moistSWE} is shown to be a particular case within a more general formulation of the moist shallow-water equations that encompasses other variants found in the literature, as well as a new formulation proposed by the authors.
\subsubsection{Governing equations}
\label{sec:mswm-eqs}
The moist shallow-water equations proposed by \citet{zerroukat2015moistSWE} may be expressed as:
\begin{align}
\label{mswm-mom}
\frac{\partial \boldsymbol{u}}{\partial t} &= -qh\boldsymbol{u}^{\perp} - \nabla B + S_u, \\
\label{mswm-h}
\frac{\partial h}{\partial t} &= -\nabla \cdot(h\boldsymbol{u}), \\
\label{mswm-temp}
\frac{\partial h \Theta}{\partial t} &= -\nabla \cdot( h \Theta\boldsymbol{u}) + hS_{\Theta}, \\
\label{mswm-tracers}
\frac{\partial h q^k}{\partial t} &= -\nabla \cdot( h q^k\boldsymbol{u}) + hS_{q^k},
\quad k =1,\ldots,3,
\end{align}
where $\boldsymbol{u}$ is the horizontal velocity, $h$ is the fluid depth, and $\Theta$ is the temperature. 
The variables $q^k$ represent moisture-related fields: $q^1$ corresponds to water vapor, $q^2$ to cloud water, and $q^3$ to rain water. The term $S_u$ is a source term for momentum, while $S_{\Theta}$ and $S_{q^k}$ are source terms for temperature and moisture, respectively.
Expressions for the source terms, describing the two-way feedback physics between vapor, cloud, and rain, can be found in \citet{zerroukat2015moistSWE} and \citet{santos2021andes}.

The potential vorticity is defined as  
\begin{equation}
\label{pot-vort}
q = \frac{\nabla \times \boldsymbol{u} + f}{h},
\end{equation}
where $f = 2\Omega \sin{\theta}$ is the Coriolis parameter, $\theta$ is the latitude, and $\Omega = 7.2921 \times 10^{-5}~\text{rad/s} $ is the Earth's rotation rate. The Bernoulli potential is given by
\begin{equation}
\label{bernoulli-pot}
B = g(h + b) + K,
\end{equation}
where $g = 9.80665~\text{m/s}^2$ is the gravitational acceleration, $b$ is the bottom topography, and $K = \frac{|\boldsymbol{u}|^2}{2}$ is the kinetic energy per unit mass.
Finally, $\boldsymbol{u}^{\perp} = \boldsymbol{k} \times \boldsymbol{u} $, where $\boldsymbol{k}$ is the local upward-pointing unit vector (i.e., normal to the sphere).

Note that Eq.\eqref{mswm-mom} represents the momentum equation of the classical shallow-water system in vector-invariant form, augmented by an additional source term. The temperature equation (Eq.\eqref{mswm-temp}) and the tracer equations (Eq.~\eqref{mswm-tracers}) are conservative advection equations with source terms, making them well-suited for finite-volume discretization. 

\subsubsection{Finite-volume moist shallow-water model}
\label{sec:mswm-discretization}
We begin the discretization of the system of equations \eqref{mswm-mom}-\eqref{mswm-tracers} using the method-of-lines approach. Following \citet{thuburn2009TRiSK, ringler2010trsk}, we first evaluate the momentum equation \eqref{mswm-ode-mom} at the edge points $\boldsymbol{x}_e$ and take its inner product with the normal vector $\boldsymbol{n}_e$. Next, we evaluate the remaining equations, \eqref{mswm-h}-\eqref{mswm-tracers}, at the cell centers $\boldsymbol{x}_i$.  
This leads to the following system of ordinary differential equations (ODEs):
\begin{align}
\label{mswm-ode-mom}
\frac{d{u}_e}{dt}(t) &= -[qh\boldsymbol{u}^{\perp} - \nabla B]_e + [S_{u}]_e, \\
\label{mswm-ode-h}
\frac{dh_i}{dt}(t) &=  -[\nabla \cdot (h\boldsymbol{u})]_i, \\
\label{mswm-ode-temp}
\frac{d (h \Theta)_i}{dt}(t)  &= -[\nabla \cdot( h \Theta\boldsymbol{u})]_i + h_i[S_{\Theta}]_i, \\
\label{mswm-ode-tracers}
\frac{d (h q^k)_i}{dt}(t) &= -[\nabla \cdot( h q^k\boldsymbol{u})]_i + h_i[S_{q^k}]_i,
\quad k = 1,\ldots,3,
\end{align}
where $[G]_k = G(\boldsymbol{x}_k, t)$ for any scalar quantity $G$, with $k = i$ (cell center) or $k = e$ (edge), and $F_e = \boldsymbol{F}(\boldsymbol{x}_e) \cdot \boldsymbol{n}_e$ for any vector quantity $\boldsymbol{F}$, where $\cdot$ denotes the standard inner product, following the notation of \citet{santos2021andes}.

We then define the full state vector as
\begin{align}
 X = \big(
 \left\{{u}_e\right\}_{e=1}^{N_e},\
 \left\{h_i\right\}_{i=1}^{N_c},\ 
 \left\{(h\Theta)_i\right\}_{i=1}^{N_c},\
 \left\{(hq^1)_i\right\}_{i=1}^{N_c}, \
 \left\{(hq^2)_i\right\}_{i=1}^{N_c}, \
 \left\{(hq^3)_i\right\}_{i=1}^{N_c} 
 \big),
 \end{align}
and express the system compactly as
 \begin{equation}
    \frac{dX}{dt}(t) = F(X(t)),
\end{equation}
where $F(X(t))$ denotes the discrete operator representing all spatial terms in the system.  
This ODE is solved using the three-stage Runge–Kutta scheme described in Appendix A.

To complete the model discretization, we define the discrete spatial operators composing $F(X(t))$. The momentum and fluid depth equations are discretized with the TRiSK scheme \citep{thuburn2009TRiSK,ringler2010trsk}, as in \citet{santos2021andes}. The temperature equation \eqref{mswm-ode-temp} and the moisture-tracer equations \eqref{mswm-ode-tracers} are discretized with either the OG or SG high-order advection schemes considered here. In both cases, following \citet{skamarockgassmann2011}, a flux-corrected transport scheme \citep{zalesak1979fct} is applied to prevent nonphysical negative values, as described in Appendix B.

TRiSK is a C-grid finite-volume and finite-difference scheme that uses the normal component of velocity at edge points, where the primal and dual grids intersect, as the prognostic wind velocity variable.  This is exactly the velocity information required by the SG scheme. However, the OG scheme requires velocity values at one or two quadrature points that do not coincide with the edge locations used by TRiSK. In this work, we adopt a least-squares reconstruction method, following \citet{peixoto2014vecrecon}, to recover the full velocity vector at the quadrature points using only the normal velocities at the edge points, achieving second-order accuracy. More details about the wind reconstruction are provided in Appendix C.

\citet{santos2021andes} applied TRiSK to the moist shallow-water model on a locally refined SCVT grid, with high resolution over the Andes, intermediate resolution over South America, and coarser global resolution, to study cloud and rain formation. They found that even smooth grid refinement could generate spurious cloud and rain. However, their study used the SG2 scheme for temperature and moisture advection, together with a monotonic flux limiter that broke total mass conservation. Here, we instead use higher-order SG and OG schemes with the more robust flux limiter described here to assess their sensitivity in simulations of cloud and rain formation over the Andes on the topography-based locally refined SCVT grids of \citet{santos2021andes}.

The study of \citet{santos2021andes} was motivated by earlier works showing that TRiSK operators may lack consistency and proper convergence on Voronoi meshes with complex cell geometries \citep{peixoto2013gridimprinting,peixoto2016consitency}. TRiSK extends the \citet{arakawa1981energyconserving} scheme from regular quadrilateral grids to arbitrary orthogonal polygonal grids while preserving key mimetic properties, such as total energy conservation and steady geostrophic modes on the $f$-sphere \citep{thuburn2009TRiSK,ringler2010trsk}. 
It was later generalized to non-orthogonal polygonal grids within a discrete exterior calculus framework \citep{thuburn2014mimeticSW,eldred2017hamiltonian,eldred2022triskdec}.

Despite its limitations, TRiSK is still used in several models, notably MPAS in its atmospheric and oceanic components \citep{skamarock2012mpas,ringler2013multiresolgrid}, as well as in other dynamical cores such as DYNAMICO \citep{dubos2015dynamico} and WAVETRiSK \citep{kevlahan2019wavetrisk}. It has also been examined in other studies \citep{weller2012compmodes,weller2012gridimprinting,weller2014NonOrthogonal}.

\section{Numerical experiments}
\label{sec:num-exp}

This section presents numerical experiments comparing the SG and proposed OG schemes. We first consider pure advection tests, following \citet{nair2010class}, and then moist shallow-water tests, following \citet{zerroukat2015moistSWE}. In the latter, the SG and OG schemes are used for tracer advection within the moist shallow-water model discretized with the TRiSK mimetic finite-volume framework.

Advection tests for the SG scheme on quasi-uniform SCVT grids were reported by \citet{skamarockgassmann2011} and are included here for comparison with the OG schemes. Since SG has not been evaluated on variable-resolution grids for pure advection, we also present results on variable-resolution SCVT grids based on the Andes topography of \citet{santos2021andes}.

For grid functions defined at control volume centers, we consider the $L_\infty$ and $L_2$ norms, as suggested in \citet{will1992swsphere}:
\begin{align*}
    \|\phi\|_{\infty} = \max_{i=1, \ldots, N_c} |\phi_i|, \quad
    \|\phi\|_{2} = \sqrt{\sum_{i=1}^{N_c} |\phi_i|^2 |\Omega_i|}.
\end{align*}
For a reference solution $\phi_{\mathrm{ref}}$, we compute the error metrics:
\begin{align*}
    E_{\infty} = \frac{\|\phi-\phi_{\mathrm{ref}}\|_{\infty}}{\|\phi_{\mathrm{ref}}\|_{\infty}}, \quad
    E_{2} = \frac{\|\phi-\phi_{\mathrm{ref}}\|_{2}}{\|\phi_{\mathrm{ref}}\|_{2}}.
\end{align*}

\subsection{Advection tests}
\label{sec:num-exp-adv}
All results for the advection equation are computed on the unit sphere.  
In this work, we denote longitude by $\lambda$ and latitude by $\theta$ on the sphere.  
The time step is chosen such that the Courant number is approximately 0.6 for all simulations.

\subsubsection{Solid-body rotation of a Gaussian hill}
\label{sec:num-exp-adv-gaussianhill}

This test considers a purely zonal wind field on the unit sphere, expressed in terms of its zonal and meridional velocity components as
\begin{equation}
u = u_{0}\cos(\theta), \quad v = 0,
\end{equation}
where \( u_{0} = \frac{2\pi}{T} \) and \( T = 5 \) time units represents the rotation period.

The experiment of \citet{will1992swsphere} describes solid-body rotation, so after one full period \(T\) the tracer field returns to its initial state. This permits a direct assessment of advection errors, and the exact solution can be computed at each time step by the method of characteristics.

\begin{figure}[!htb]
    \centering
    \includegraphics[width=1.00\linewidth]{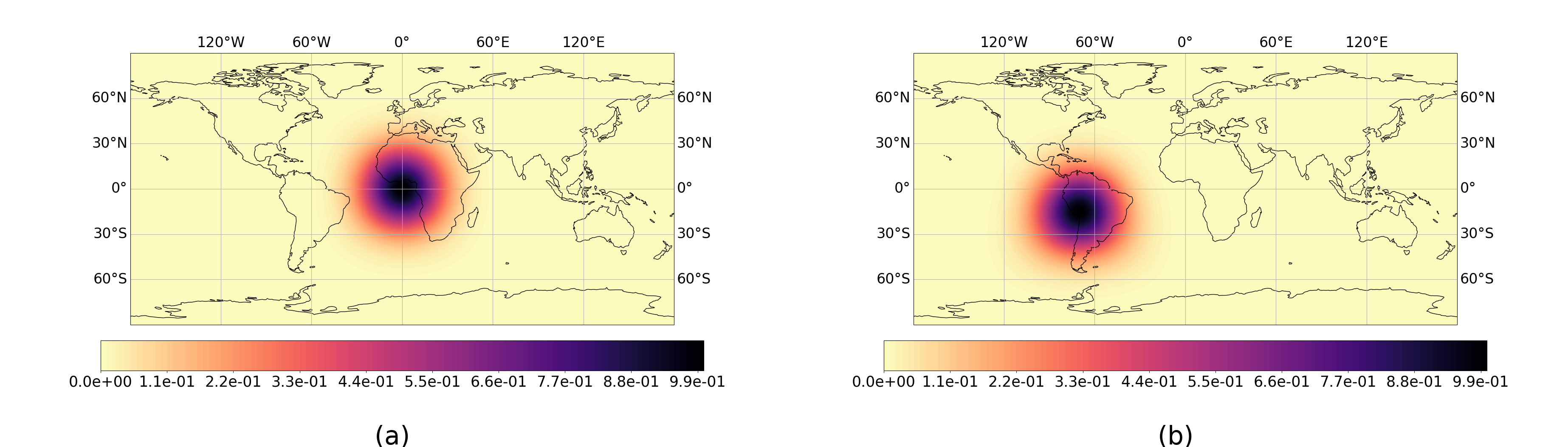}
    \caption{Initial condition for the zonal wind test case: a Gaussian hill centered at $(\lambda_0,\theta_0)=(0,0)$ (a) in the quasi-uniform grid and at $(\lambda_0,\theta_0)=\left(\tfrac{7\pi}{18}, -\tfrac{\pi}{12}\right)$ in the variable-resolution setup (b), so that the latter is located over the Andes.}
    \label{fig:adv-zonal-gaussian}
\end{figure}

The initial tracer field consists of a single Gaussian hill defined in Cartesian coordinates \((x, y, z)\) as
\begin{equation}
\phi(x, y, z) = \exp\left[-5\left((x - x_{0})^{2} + (y - y_{0})^{2} + (z - z_{0})^{2}\right)\right],
\end{equation}
centered at the point $(x_0,y_0,z_0)$ corresponding to longitude and latitude \((\lambda_{0}, \theta_{0}) = (0, 0)\).
In the uniform grid configuration, we set the center at 
$(\lambda_{0}, \theta_{0}) = (0, 0)$ (Figure \ref{fig:adv-zonal-gaussian}a), whereas for the variable-resolution 
grids we use $(\lambda_{0}, \theta_{0}) = \left(\frac{7\pi}{18}, -\frac{\pi}{12}\right)$ (Figure \ref{fig:adv-zonal-gaussian}b), 
so that the initial Gaussian hill is centered over the Andes.

\begin{figure}[!htb]
	\centering
    \includegraphics[width=1.00\linewidth]{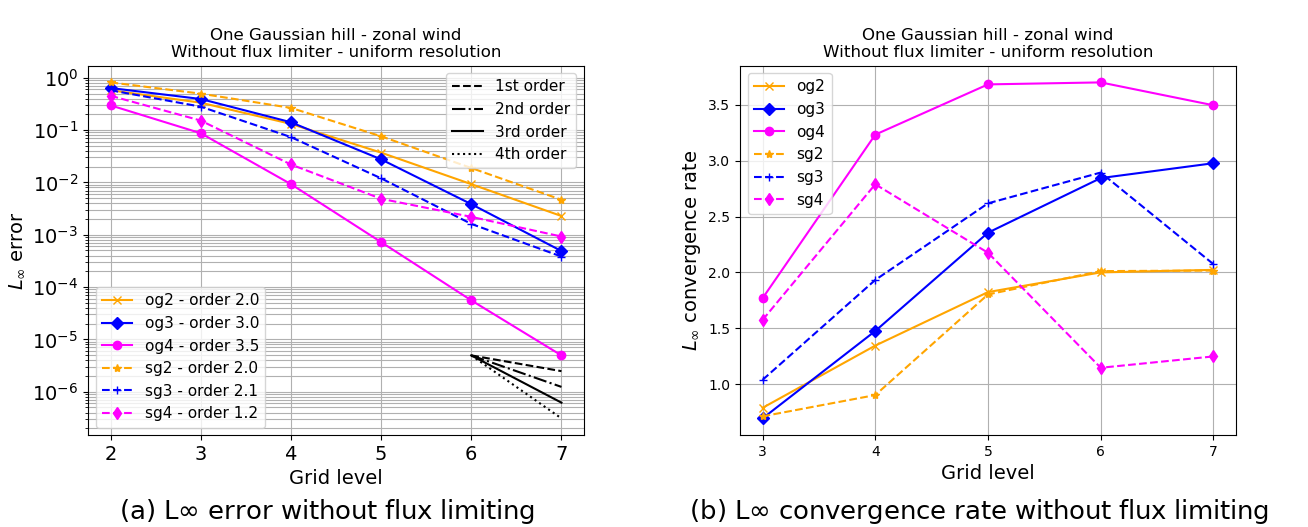}
\caption{$L_\infty$ errors (a) and convergence rates (b) for the advection of a Gaussian hill under a zonal wind, computed with OG (solid) and SG (dashed) schemes on quasi-uniform SCVT grids (levels 2-7) without flux limiting. Line colors indicate the nominal order of accuracy. The Gaussian hill is initially centered at latitude and longitude $(0,0)$ (Figure \ref{fig:adv-zonal-gaussian}a)}
\label{fig:adv-zonal-gaussian-unif-linf}
\end{figure}

We use SG3 with $\beta = 1$ in this simulation.
Figure~\ref{fig:adv-zonal-gaussian-unif-linf} presents the $L_{\infty}$ error convergence and its convergence rate on the quasi-uniform SCVT grid, without flux limiting.  
From this figure, we note that both OG2 and SG2 achieve second-order convergence, with the $L_\infty$ errors of OG2 being smaller. OG3 reaches third-order convergence; however, SG3 exhibits third-order convergence only up to level 6, after which the convergence rate decays toward second order. The $L_\infty$ errors of SG3 remain smaller than those of OG3 up to grid level 6, and at level 7 they match OG3. OG4 reaches almost 3.5th-order convergence, while SG4 achieves only first-order convergence.

\begin{figure}[!htb]
	\centering
    \includegraphics[width=1.00\linewidth]{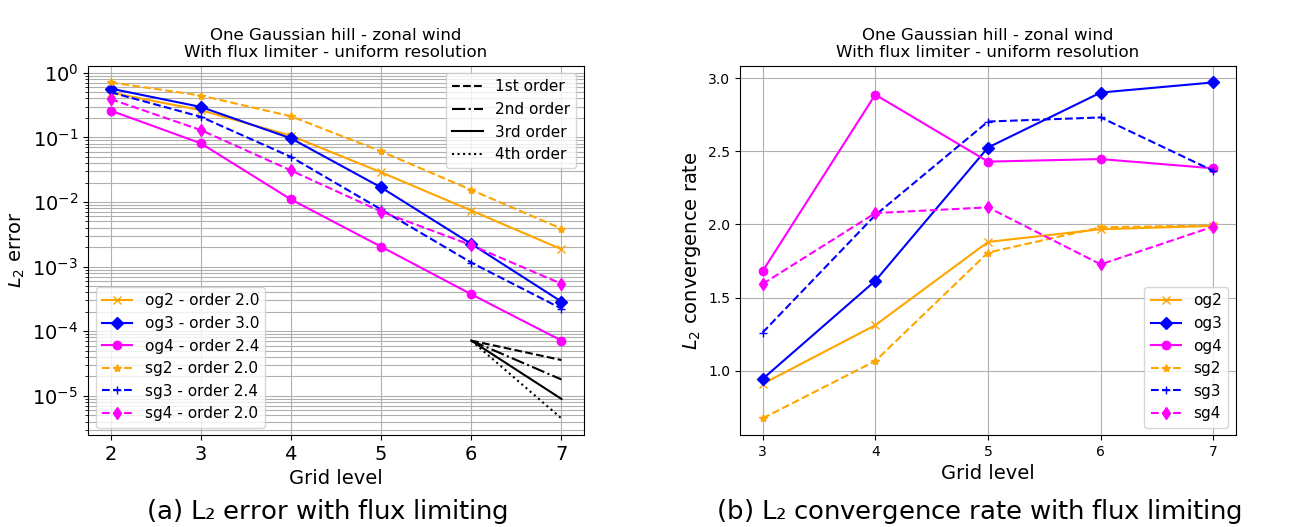}
\caption{$L_2$ errors (a) and convergence rates (b) for the advection of a Gaussian hill under a zonal wind, computed with OG (solid) and SG (dashed) schemes on quasi-uniform SCVT grids (levels 2-7) with flux limiting. Line colors indicate the nominal order of accuracy. The Gaussian hill is initially centered at latitude and longitude $(0,0)$  (Figure \ref{fig:adv-zonal-gaussian}a).}
     \label{fig:adv-zonal-gaussian-unif-l2-mono}
\end{figure}

The $L_2$ errors and convergence rates on the quasi-uniform SCVT grid with flux limiting are shown in Figure~\ref{fig:adv-zonal-gaussian-unif-l2-mono}. These results are consistent with the $L_\infty$ errors without flux limiting: OG2 yields smaller errors than SG2, with both remaining second-order accurate. In the $L_2$ norm, OG3 again achieves third-order convergence, while SG3 reaches about 2.5th order and has smaller errors than OG3 except at grid level 7, where they become comparable. SG4 performs better with flux limiting, attaining nearly second-order convergence, but OG4 still yields smaller errors, with a convergence rate of about 2.5.

Overall, these results indicate that OG2 generally has smaller errors than SG2, while both remain second-order accurate. The only SG scheme that outperforms its OG counterpart is SG3, although its convergence order decreases at finer grids, where the errors become comparable. The most problematic SG scheme is SG4, which performs poorly without flux limiting. By contrast, OG4 consistently yields the smallest errors overall, both with and without flux limiting.

\begin{figure}[!htb]
	\centering
    \includegraphics[width=1.00\linewidth]{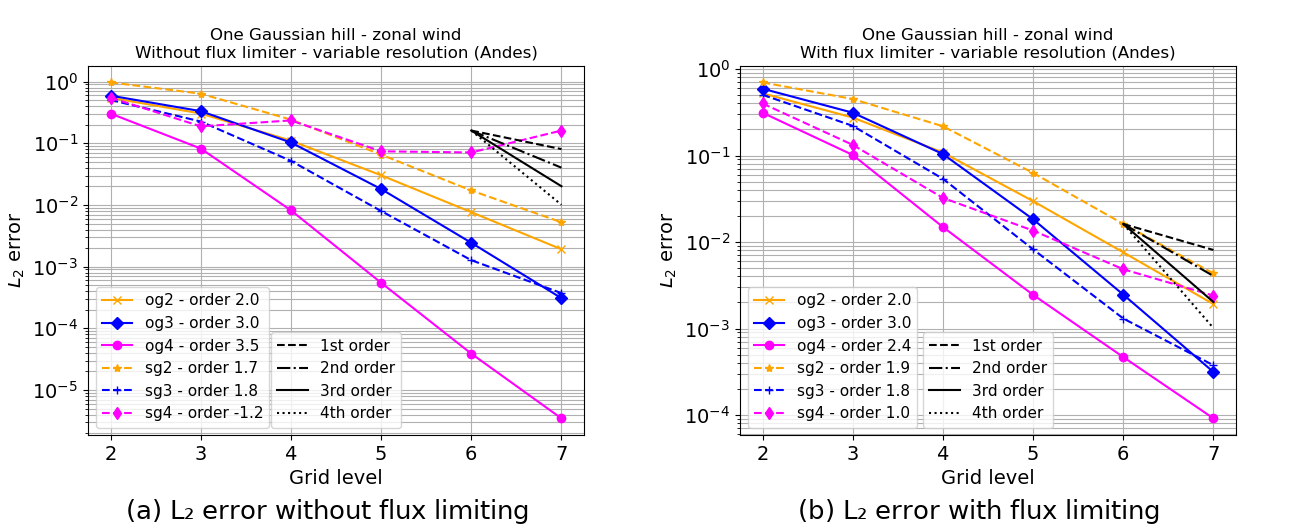}
\caption{$L_2$ error convergence without (a) and with flux limiting (b) for the advection of a Gaussian hill under a zonal wind, computed with OG (solid) and SG (dashed) schemes on variable-resolution SCVT grids (levels 2-7). Line colors indicate the nominal order of accuracy. The Gaussian hill is initially centered at latitude and longitude $\left(\frac{7\pi}{18}, -\frac{\pi}{12}\right)$, located over the Andes  (Figure \ref{fig:adv-zonal-gaussian}b).}
\label{fig:adv-zonal-gaussian-ref-l2}
\end{figure}

Next, we repeat the test on a variable-resolution grid refined according to the Andes topography (Figure~\ref{fig:scvt_grid}b). The $L_2$ error convergence without and with flux limiting is shown in Figure~\ref{fig:adv-zonal-gaussian-ref-l2}. OG4 is particularly noteworthy, as it preserves the best overall accuracy on the variable-resolution grid, showing strong robustness even under topography-based local refinement. With or without flux limiting, the nominally second- and third-order schemes give results similar to those on quasi-uniform SCVT grids, indicating low sensitivity to the more challenging topography-based refinement. In contrast, SG4 produces much larger errors, failing to converge without flux limiting and achieving only first-order convergence with it.

\begin{figure}[!htb]
	\centering
    \includegraphics[width=1.\linewidth]{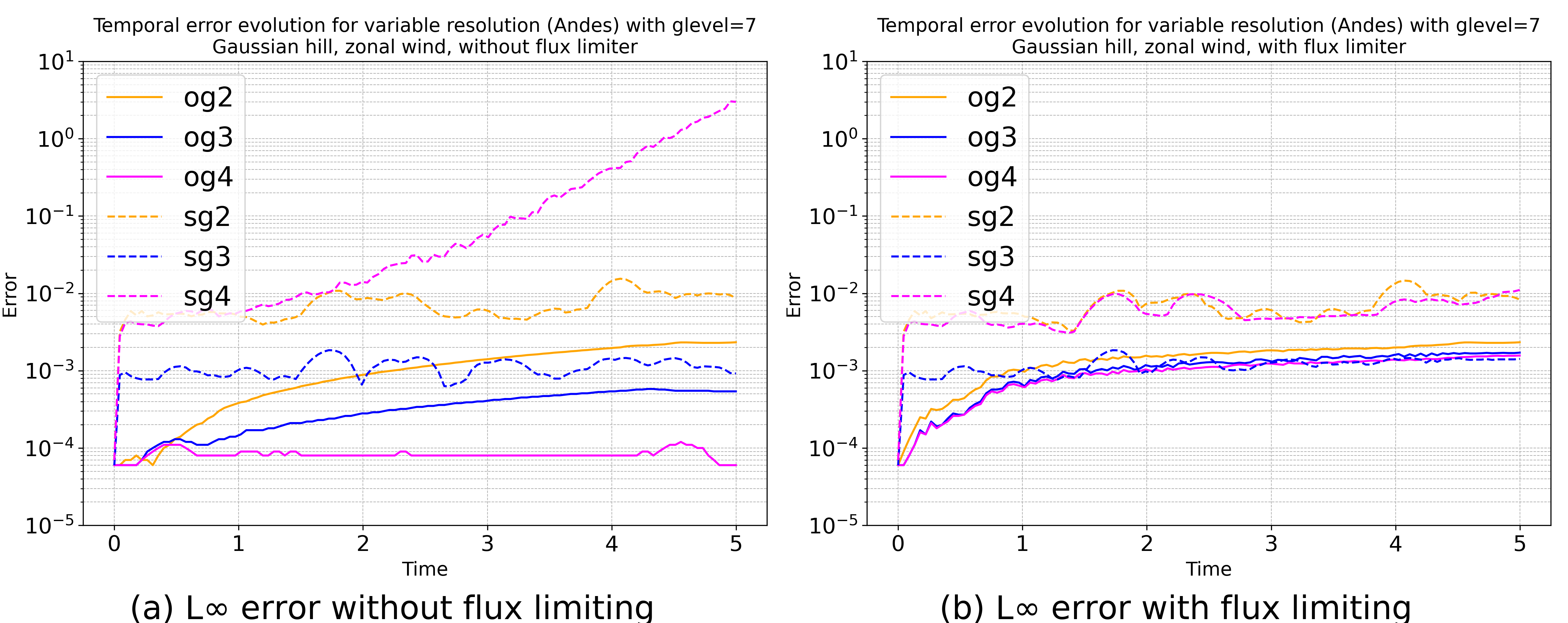}
\caption{Time evolution of the $L_\infty$ error for the advection of a Gaussian hill under a zonal wind over 5 time units, without (a) and with (b) flux limiting, computed using OG (solid) and SG (dashed) schemes on level 7 variable-resolution SCVT grid. Line colors indicate the nominal order of accuracy. The Gaussian hill is initially centered at latitude and longitude $\left(\frac{7\pi}{18}, -\frac{\pi}{12}\right)$, located over the Andes  (Figure \ref{fig:adv-zonal-gaussian}b).}
\label{fig:adv-zonal-gaussian-ref-linf}
\end{figure}

To investigate further, we show the time evolution of the $L_{\infty}$ error on the variable-resolution SCVT grid at level 7, without and with flux limiting, in Figure~\ref{fig:adv-zonal-gaussian-ref-linf}. At this level, OG3 has smaller errors than SG3 without flux limiting, while their errors are very similar with flux limiting. As expected, SG2 has larger errors than OG2 in both cases.

The $L_{\infty}$ errors of SG4 without flux limiting increase over time, indicating possible instability. With flux limiting, the errors remain stable and SG4 behaves similarly to SG2. Note that SG2 and SG4 are the only schemes without upwind bias, making them more prone to numerical instabilities. The added numerical diffusion from flux limiting explains the stable $L_{\infty}$ errors of SG4 in this case.

\subsubsection{Deformational flow with two Gaussian hills}
\label{sec:num-exp-adv-twohills}

This test follows the deformational flow experiment proposed by \citet{nair2010class}, which considers two Gaussian hills as the initial condition:
\begin{equation}
\phi(\lambda,\theta) = h_{1}(\lambda,\theta) + h_{2}(\lambda,\theta),
\end{equation}
where each Gaussian hill $h_i$ is defined in Cartesian coordinates as
\begin{equation}
h_{i} = \exp\left[-5\left((x - x_{i})^{2} + (y - y_{i})^{2} + (z - z_{i})^{2}\right)\right].
\end{equation}

\begin{figure}[!htb]
    \centering
    \includegraphics[width=1.00\linewidth]{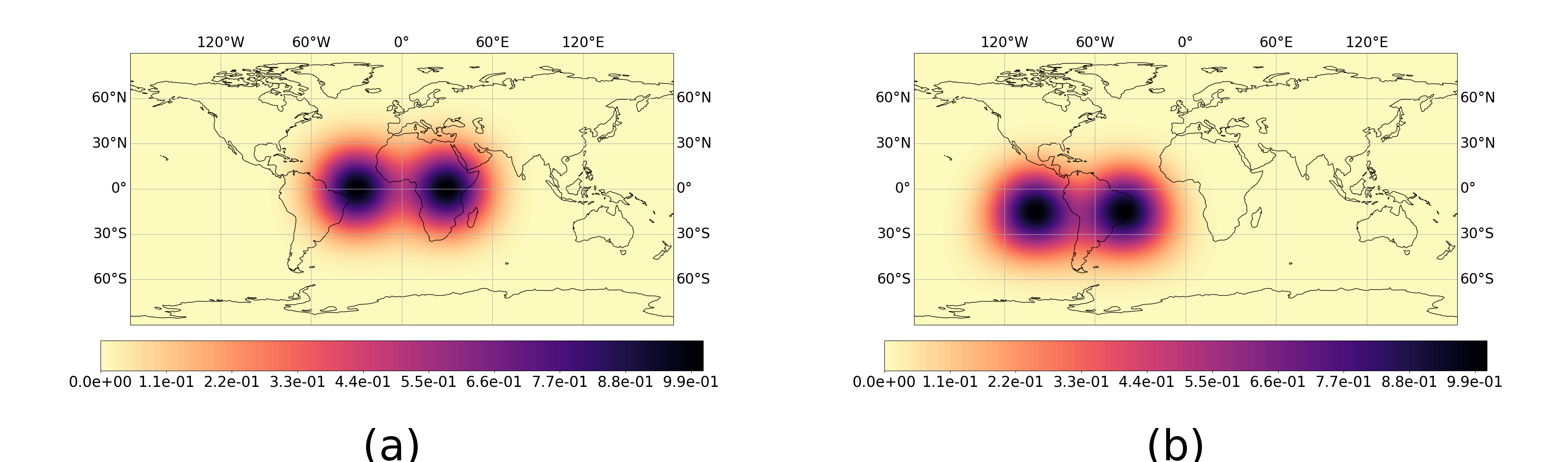}
    \caption{Initial condition for the deformational flow test case: two Gaussian hills on (a) the quasi-uniform SCVT grid and (b) the variable-resolution SCVT grid. In the variable-resolution setup, the hills are positioned over the Andes.}    \label{fig:adv-deform-gaussian}
\end{figure}

For the quasi-uniform SCVT grid, the Gaussian hills are centered at $(\lambda_{1}, \theta_{1}) = (-\pi/6, 0)$ and $(\lambda_{2}, \theta_{2}) = (\pi/6, 0)$, respectively. For the variable-resolution grid (Figure~\ref{fig:adv-deform-gaussian}a), we adopt $(\lambda_{1}, \theta_{1}) = (-\frac{5\pi}{9}, -\frac{\pi}{12})$ and $(\lambda_{2}, \theta_{2}) = (-\frac{2\pi}{9},  -\frac{\pi}{12})$, and the Gaussian hills are symmetrically distributed over the Andes in the longitudinal direction (Figure \ref{fig:adv-deform-gaussian}b).

The velocity field on the unit sphere is given by
\begin{align}
\label{eq:deform-u}
\displaystyle u(\lambda,\theta,t) &= k \sin^{2}(\lambda' + \pi) \sin(2\theta) \cos(\pi t/T) + \frac{2 \pi \cos(\theta)}{T},\\
\label{eq:deform-v}
\displaystyle v(\lambda,\theta,t) &= k \sin\left(2(\lambda' + \pi)\right) \cos(\theta) \cos(\pi t/T),
\end{align}
where $\lambda' = \lambda - 2 \pi/T$, and the additional zonal wind term is included to avoid error cancellation \citep{nair2010class}.  
This velocity field deforms and transports the Gaussian hills, which return to their initial positions after one full period $T$, allowing for the computation of numerical errors in the advection schemes. We adopt $T = 5$ time units, and SG3 uses $\beta = 1$.

\begin{figure}[!htb]
	\centering
	\includegraphics[width=1.00\linewidth]{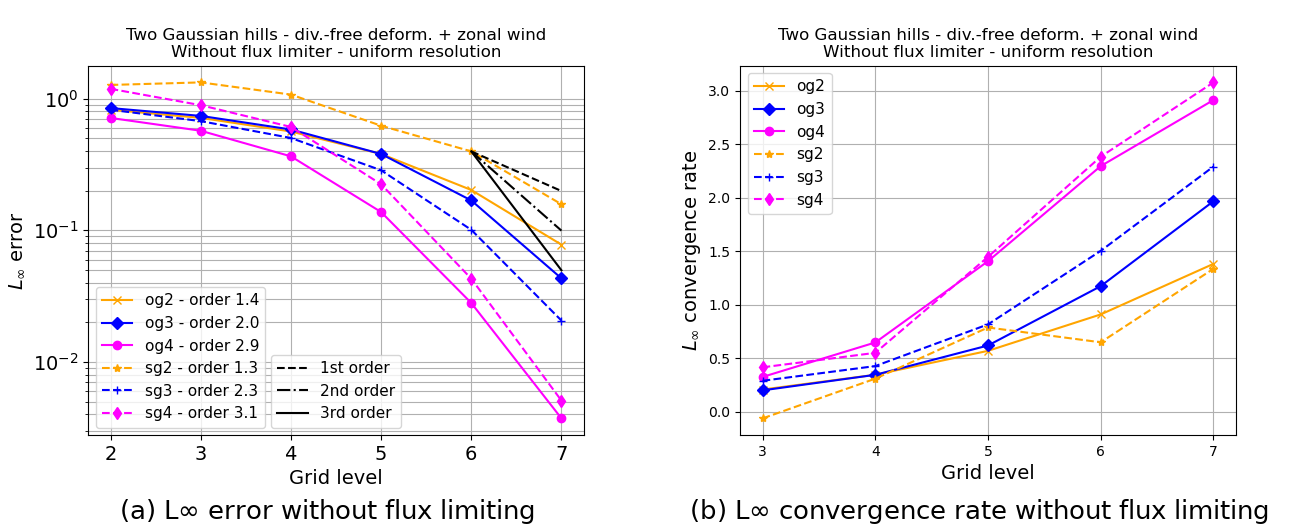}
\caption{$L_\infty$ errors (a) and convergence rates (b) for the advection of two Gaussian hills under a deformational flow velocity field, computed with OG (solid) and SG (dashed) schemes on quasi-uniform SCVT grids (levels 2-7) without flux limiting. Line colors indicate the nominal order of accuracy. The initial condition is shown in Figure \ref{fig:adv-deform-gaussian}a.}
\label{fig:adv-deform-gaussian-unif-linf}
\end{figure}

\begin{figure}[!htb]
	\centering
    	\includegraphics[width=1.00\linewidth]{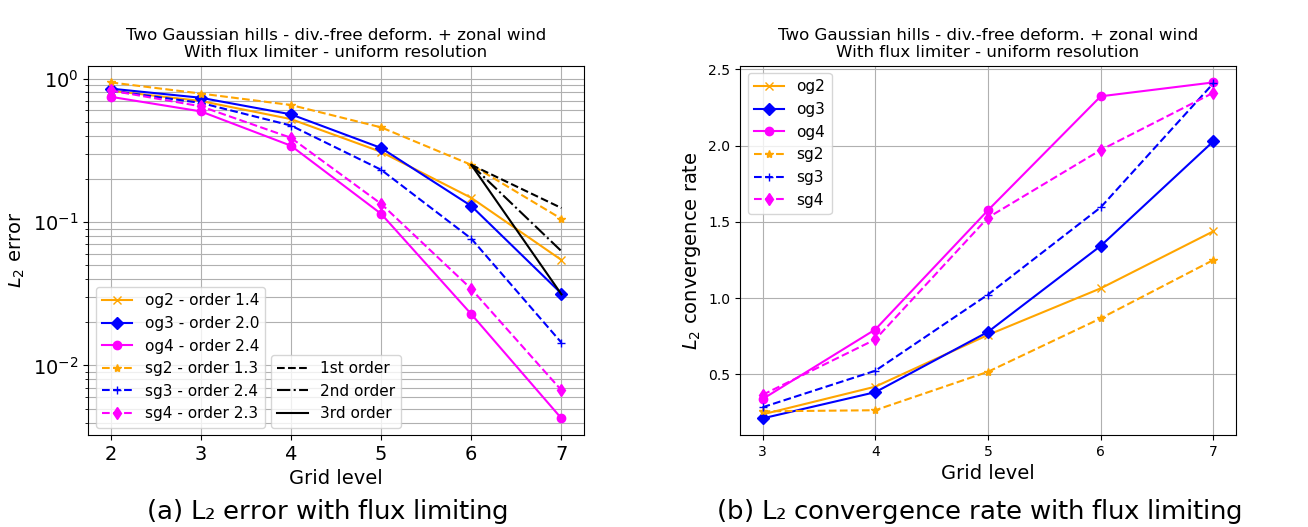}
\caption{$L_2$ errors (a) and convergence rates (b) for the advection of two Gaussian hills under a deformational flow velocity field, computed with OG (solid) and SG (dashed) schemes on variable-resolution SCVT grids (levels 2-7) with flux limiting. Line colors indicate the nominal order of accuracy.  The initial condition is shown in Figure \ref{fig:adv-deform-gaussian}a.}
     \label{fig:adv-deform-gaussian-unif-l2}
\end{figure}

Figure~\ref{fig:adv-deform-gaussian-unif-linf} shows the $L_\infty$ errors without flux limiting and the corresponding convergence rates. OG2 again yields smaller errors than SG2, whereas SG3 yields smaller errors than OG3. Unlike in the previous test case, SG4 performs comparably to OG4, with OG4 showing only slightly smaller errors. Notably, SG4 had previously shown loss of convergence in the advection of a single Gaussian hill under zonal wind, especially without flux limiting, but this behavior is not observed here. Figure~\ref{fig:adv-deform-gaussian-unif-l2} shows the $L_2$ errors with flux limiting, which exhibit similar qualitative behavior.

Regarding the convergence rates, Figure~\ref{fig:adv-deform-gaussian-unif-linf}a shows some improvement, although they remain below their nominal values. The nominally fourth-order schemes, OG4 and SG4, approach third-order convergence, while the nominally third-order schemes, SG3 and OG3, converge at about second order. The nominally second-order schemes, OG2 and SG2, reach convergence rates near 1.4.

\begin{figure}[!htb]
	\centering
    \includegraphics[width=1.00\linewidth]{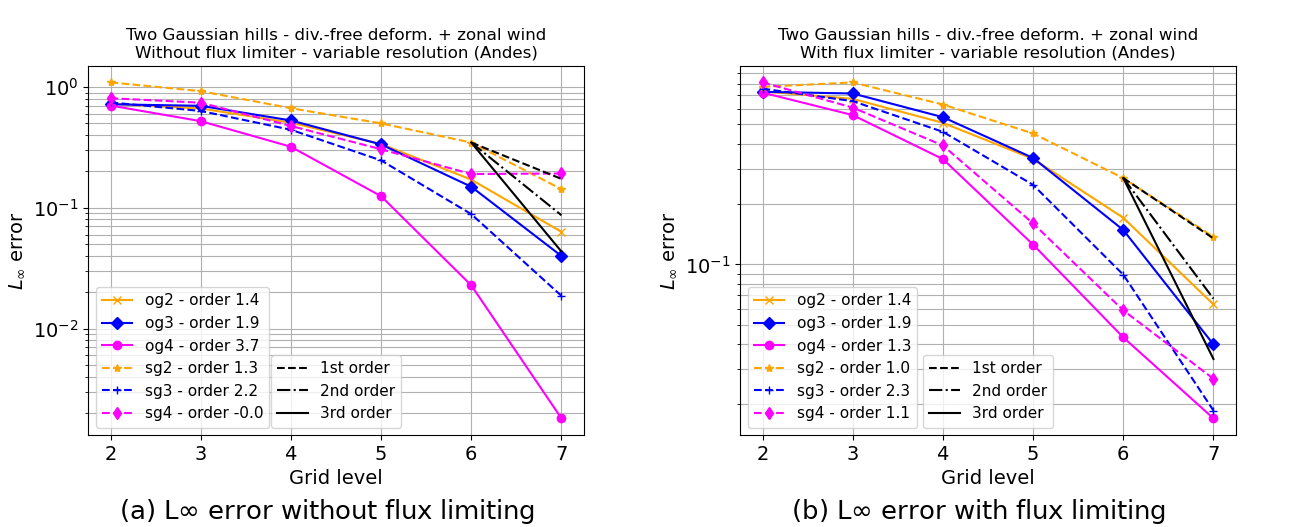}
\caption{$L_{\infty}$ error convergence without (a) and with flux limiting (b) for the advection of a Gaussian hill under a deformational flow velocity field, computed with OG (solid) and SG (dashed) schemes on variable-resolution SCVT grids (levels 2-7). Line colors indicate the nominal order of accuracy. The initial condition is shown in Figure \ref{fig:adv-deform-gaussian}b.}
\label{fig:adv-deform-gaussian-ref-linf}
\end{figure}

We again use variable-resolution grids to assess the sensitivity of each scheme to local refinement over the Andes. Figure~\ref{fig:adv-deform-gaussian-ref-linf} shows the $L_\infty$ errors without and with flux limiting.
In these cases, the OG schemes show little sensitivity to the variable-resolution grid.
Without flux limiting, OG4 reaches a convergence rate close to 3.7, while with flux limiting it drops to about 1.3, as expected from the effect of flux limiting.
Even so, OG4 yields the smallest errors overall, except at grid level 7, where its errors match those of SG3. Likewise, OG2 consistently produces smaller errors than SG2, with and without flux limiting. SG3 also performs very well with flux limiting, surpassing OG3 and matching OG4 at the highest resolution.

The only convergence issue arises with SG4, which does not converge without flux limiting. As in the previous zonal-wind test with a single Gaussian hill, this is likely due to the lack of upwind bias combined with the larger stencil.

\subsubsection{Deformational flow with two slotted cylinders}
\label{sec:num-exp-adv-slottedcylinders}

\begin{figure}[!ht]
	\centering
    \includegraphics[width=1\linewidth]{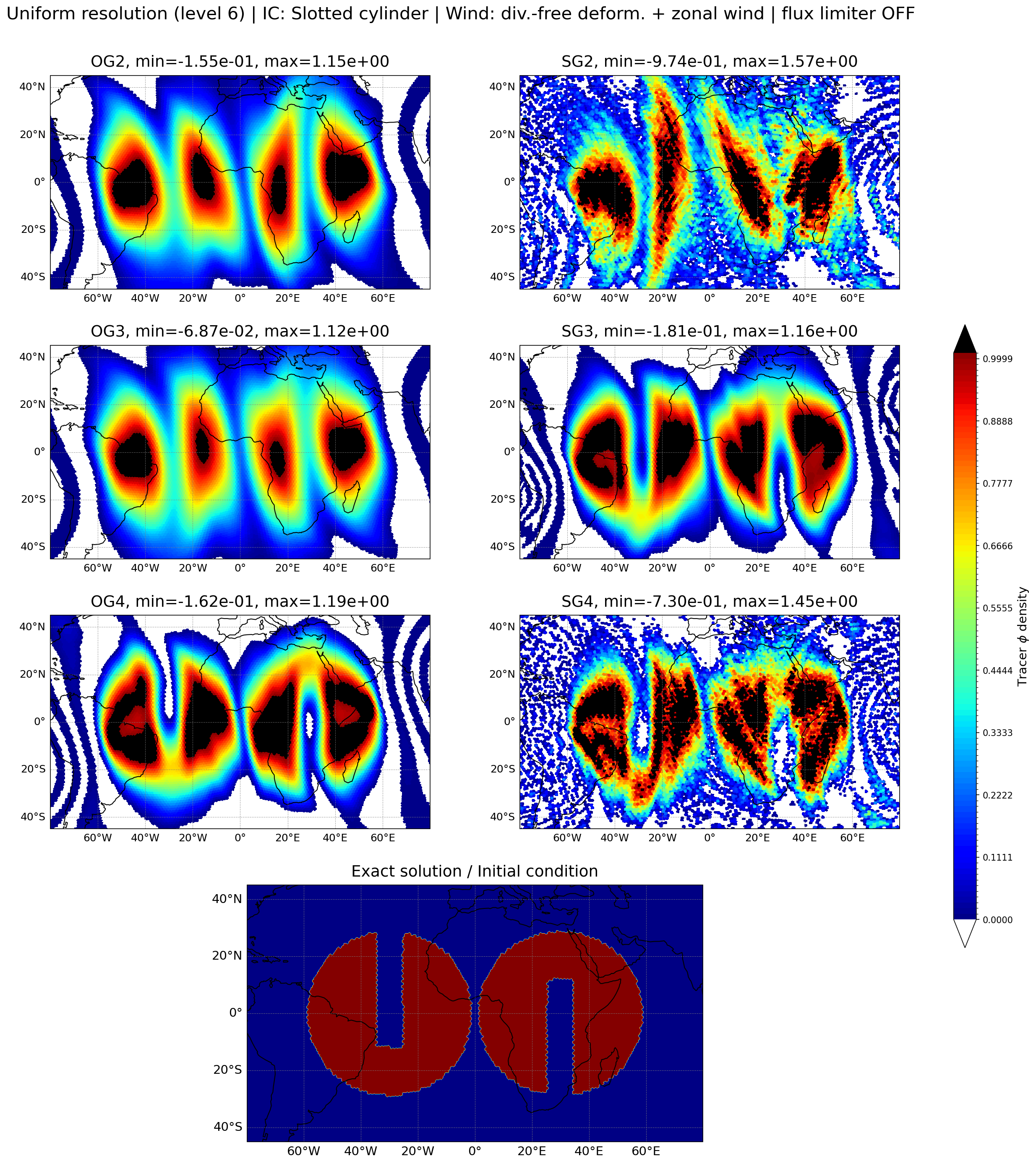}
\caption{Slotted cylinder advection under the deformational flow velocity field on a quasi-uniform SCVT grid (level~6) without flux limiting. OG schemes are shown on the left, SG schemes on the right, and the exact solution (initial condition) is shown in the bottom center panel.}
\label{fig:adv-deform-cylinders-unif-nolim}
\end{figure}

\begin{figure}[!ht]
	\centering
    \includegraphics[width=1\linewidth]{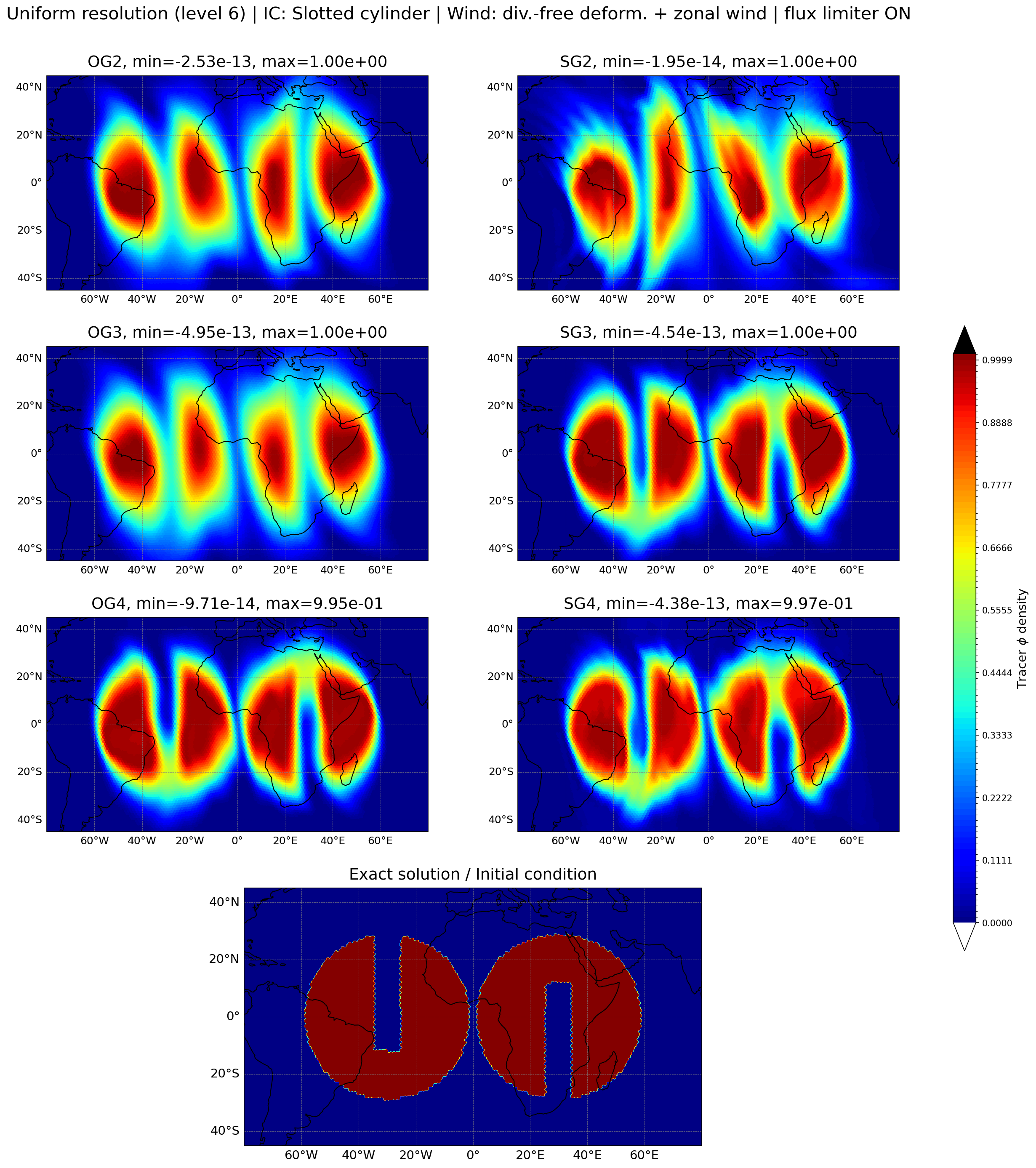}
\caption{Slotted cylinder advection under the deformational flow velocity field on a quasi-uniform SCVT grid (level~6) with flux limiting. OG schemes are shown on the left, SG schemes on the right, and the exact solution (initial condition) is shown in the bottom center panel.}
\label{fig:adv-deform-cylinders-unif-lim}
\end{figure}

\begin{figure}[!ht]
	\centering
    \includegraphics[width=1\linewidth]{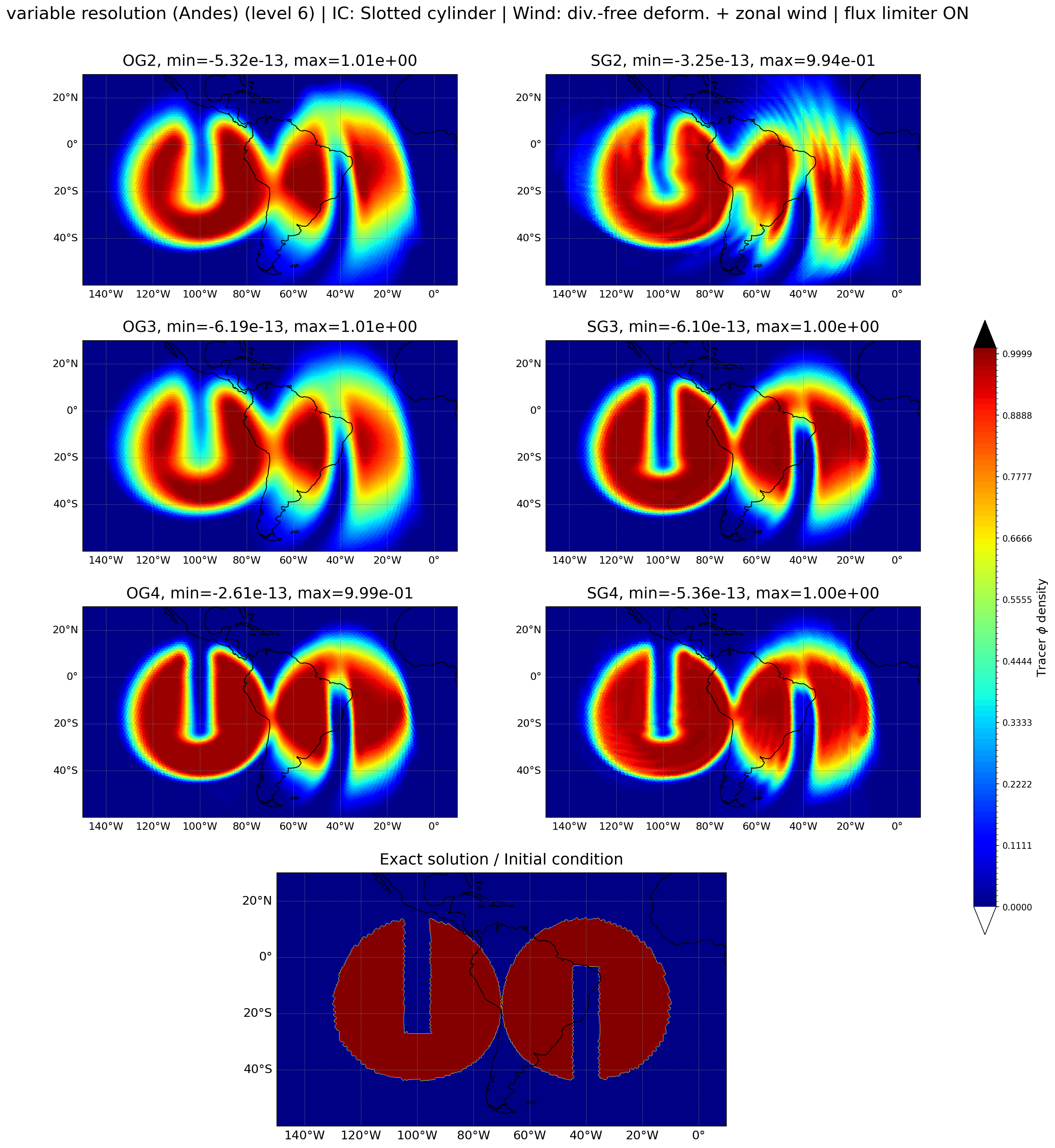}
\caption{Slotted cylinder advection under the deformational flow velocity field on a variable-resolution SCVT grid (level~6) with flux limiting. OG schemes are shown on the left, SG schemes on the right, and the exact solution (initial condition) appears in the bottom center panel.}
\label{fig:adv-deform-cylinders-var-lim}
\end{figure}

This test uses the same deformational flow as in the previous experiment (Eqs.~\ref{eq:deform-u} and \ref{eq:deform-v}), with period $T = 5$, but the initial condition now consists of two slotted cylinders. Their analytical expressions are given by \citet{nair2010class}, and the initial profile is shown in the bottom panel of Figure~\ref{fig:adv-deform-cylinders-unif-nolim}. For the variable-resolution setup, the initial condition is shifted to be symmetric over the Andes, as shown in the bottom panel of Figure~\ref{fig:adv-deform-cylinders-var-lim}.

Because the slotted-cylinder field is discontinuous, this test is well suited for evaluating flux limiters. The discontinuities generate numerical oscillations, allowing us to assess how well each scheme preserves monotonicity and suppresses spurious extrema. 
For this simulation, SG3 uses $\beta = 0.25$.

Figure~\ref{fig:adv-deform-cylinders-unif-nolim} shows the final cylinder profiles for all OG and SG schemes on a uniform level-6 grid without flux limiting. All schemes produce new extrema, with the largest oscillations in the centered schemes SG2 and SG4. With flux limiting (Figure~\ref{fig:adv-deform-cylinders-unif-lim}), these oscillations are suppressed and no new extrema appear. Among the limited schemes, OG2 outperforms SG2, SG3 outperforms OG3, and OG4 outperforms SG4, consistent with the previous tests.
Figure~\ref{fig:adv-deform-cylinders-var-lim} shows the variable-resolution results at refinement level 6 with flux limiting. The same qualitative behavior is observed: OG4 performs best overall, followed by SG3. Small oscillations remain visible in the centered schemes SG2 and SG4.

\subsection{Moist shallow-water tests}
\label{sec:num-exp-mswm}
In this subsection, we present numerical experiments for the moist shallow-water model of \citet{zerroukat2015moistSWE} to assess the higher-order advection schemes within the TRiSK-based mimetic finite-volume framework.

All simulations are performed on the sphere with Earth’s radius $R = 6.371 \times 10^6\ \text{m}$. The SG and OG schemes are used to advect moisture tracers and temperature, while the fluid depth and normal velocity equations are solved with TRiSK \citep{thuburn2009TRiSK,ringler2010trsk}.
Unlike SG, the OG scheme also requires the wind reconstruction of Appendix C to obtain velocities at quadrature points. The SG3 scheme uses $\beta = 0.25$ in the simulations presented here.

On variable-resolution grids, we use fourth-order hyperdiffusion with a coefficient based on the local cell diameter, following \citet{santos2021andes}, to ensure numerical stability. No hyperdiffusion is applied on quasi-uniform SCVT grids unless otherwise mentioned.

\subsubsection{Moist steady geostrophic flow}
\label{sec:num-exp-mswm-tc2}
The first moist shallow-water test case is the moist steady geostrophic flow of \citet{zerroukat2015moistSWE}, which extends the steady geostrophic flow test of \citet{will1992swsphere} to the moist shallow-water system. The initial conditions are in exact balance and remain steady, making this case suitable for assessing numerical errors; ideally, no rain or cloud should form. In practice, the simulations produced small cloud-related errors but no rain.

\begin{figure}[!htb]
    \centering
    \includegraphics[width=\linewidth]{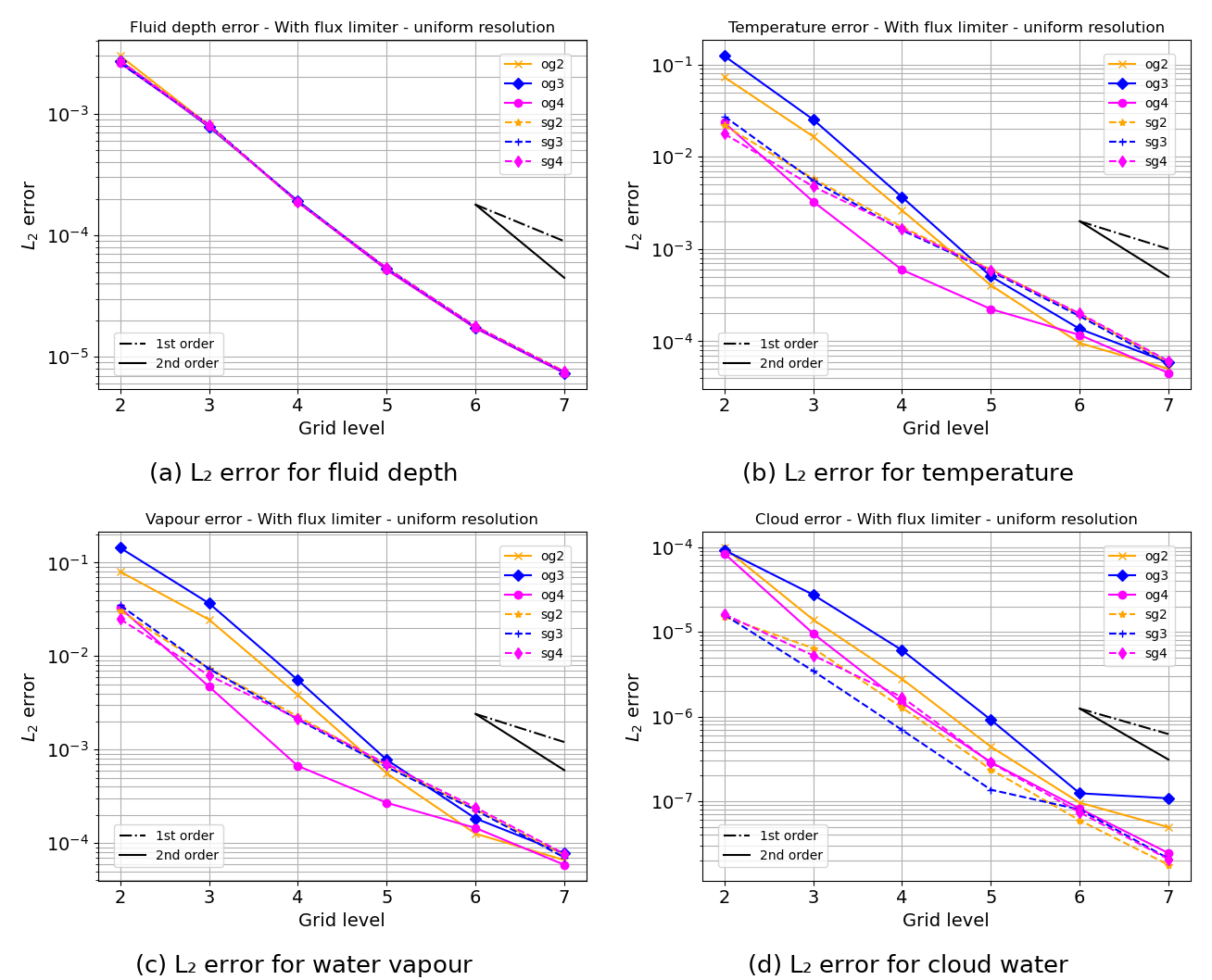}
    \caption{$L_{2}$ error convergence for the moist shallow-water model in the moist steady geostrophic test case from \citet{zerroukat2015moistSWE}, considering 12 days of simulation.
    Results are shown for simulations with flux limiting applied to the fluid depth field, using the OG (solid lines) and SG (dashed lines) schemes on quasi-uniform SCVT grids (levels 2-7). 
    Line colors indicate the nominal order of accuracy.}
    \label{fig:mswm-tc2}
\end{figure}

\begin{figure}[!htb]
    \centering
        \includegraphics[width=\linewidth]{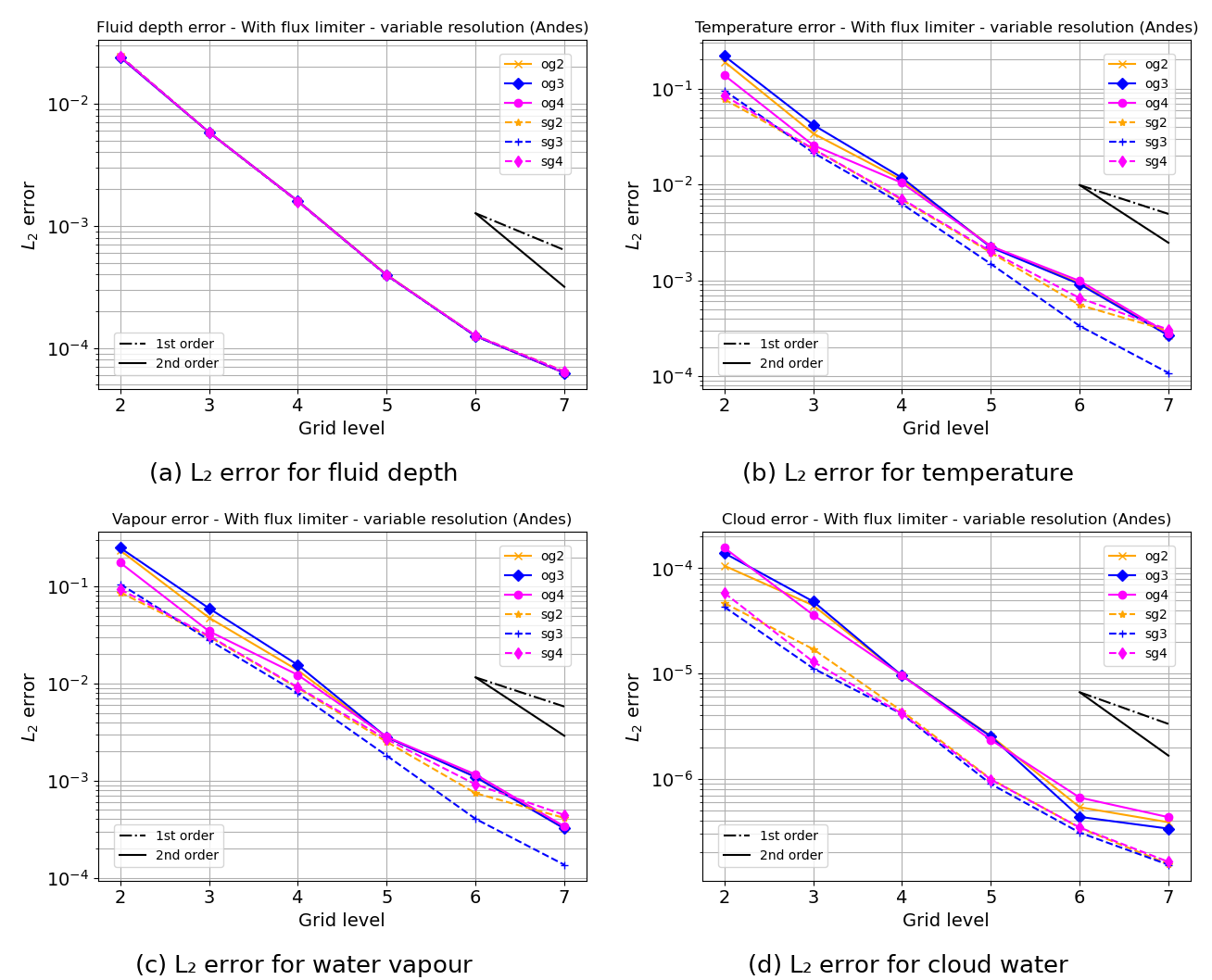}
    \caption{$L_{2}$ error convergence for the moist shallow-water model in the moist steady geostrophic test case from \citet{zerroukat2015moistSWE}, considering 12 days of simulation.
    Results are shown for simulations with flux limiting applied to the fluid depth field, using the OG (solid lines) and SG (dashed lines) schemes on variable-resolution SCVT grids (levels 2-7). 
    Line colors indicate the nominal order of accuracy.}
    \label{fig:mswm-tc2-ref}
\end{figure}

Figure~\ref{fig:mswm-tc2} shows the $L_2$ errors in fluid depth, temperature, water vapor, and cloud water after 12 days on uniform-resolution SCVT grids with flux limiting, using both OG and SG schemes. Figure~\ref{fig:mswm-tc2-ref} shows the corresponding results on variable-resolution grids. In both cases, grid levels 2 to 7 are considered. The time step at level 2 is $5120$ seconds for the uniform grid and $2560$ seconds for the variable-resolution grid, and is halved at each finer level.

From Figure~\ref{fig:mswm-tc2}, we observe that on quasi-uniform SCVT grids the $L_2$ errors of fluid depth change little across advection schemes and initially show second-order convergence, later degrading toward first order. This is expected, since the TRiSK divergence operator is inconsistent and does not fully converge on SCVT grids \citep{peixoto2016consitency}. In fact, the corresponding $L_\infty$ error of fluid depth (not shown) exhibits zero-order convergence. Because the tracer equations are written in flux form and depend on fluid depth, the tracer fields are also affected by this limitation, even with higher-order advection schemes. Accordingly, temperature and water vapor errors initially show second-order or higher convergence, but later decrease to between first and second order. Cloud errors show a similar trend, although SG3 loses convergence from levels 5 to 6 and recovers it from 6 to 7, while OG3 loses convergence from levels 6 to 7.

In general, from Figure~\ref{fig:mswm-tc2} we can see that the $L_2$ errors for the OG schemes are generally larger than those for the SG schemes for temperature, water vapor, and cloud water, respectively. We attribute this difference to the additional errors introduced by the wind reconstruction procedure required by the OG schemes.
Among the OG schemes, OG4 yields the most accurate results in the temperature and water vapor fields (Figures~\ref{fig:mswm-tc2}b and \ref{fig:mswm-tc2}c), particularly at coarser resolutions, where it is more accurate than the SG schemes.
For the cloud field, the SG schemes generally produce more accurate results, especially at coarse resolutions. At the highest resolutions, however, OG4 yields results that are comparable to those of the SG schemes.
Overall, the SG schemes exhibit very similar error magnitudes for temperature and water vapor. Among them, SG3 is more accurate than the other SG variants, particularly for the cloud field at coarser resolutions.

\begin{figure}[!htb]
	\centering
    \includegraphics[width=1.\linewidth]{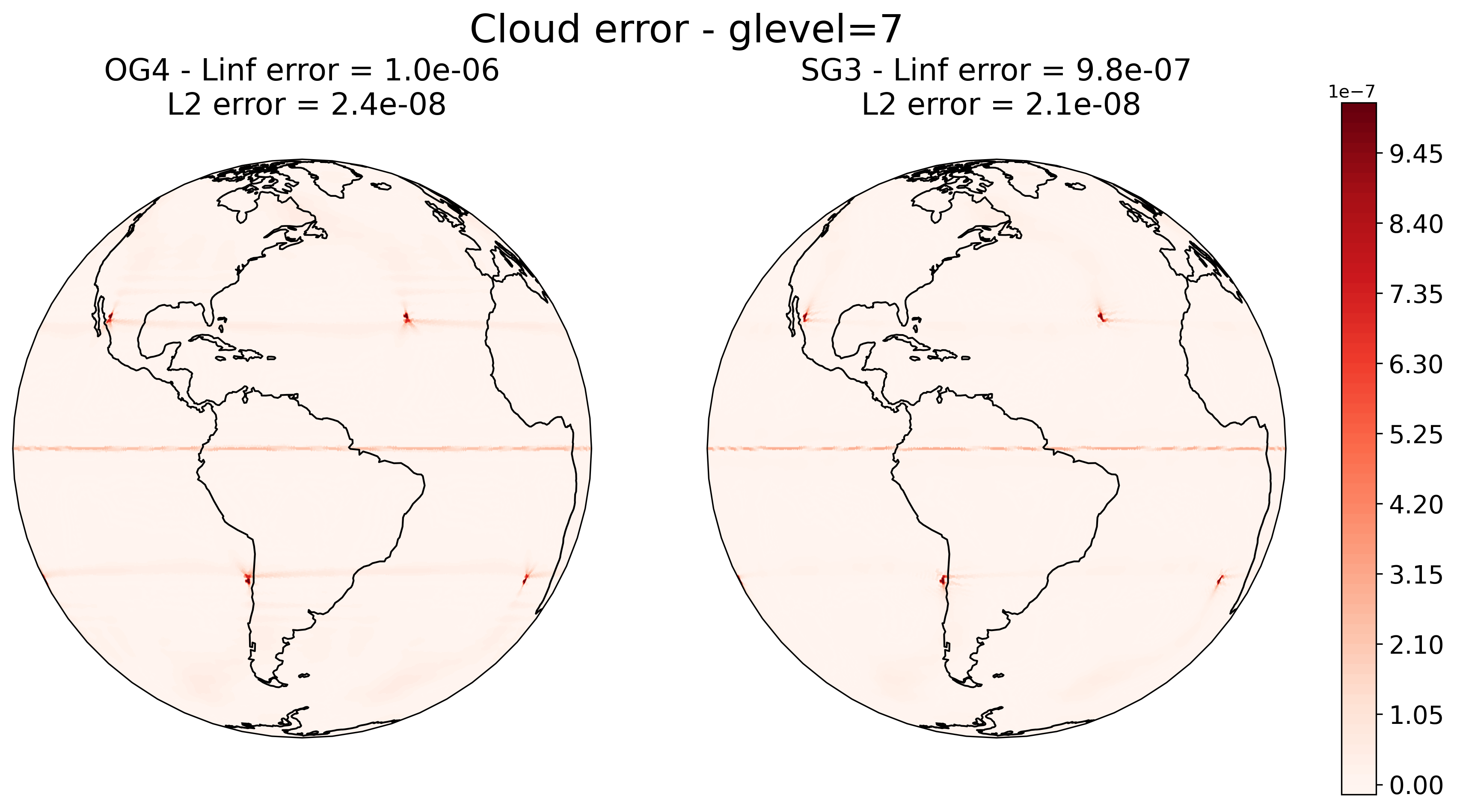}
\caption{Spatial distribution of the errors for the cloud water field at the uniform SCVT grid level 7, comparing the OG4 (right) and SG3 (left) schemes, for the moist steady geostrophic flow in the moist shallow water model.}
\label{fig:mswm-tc2-cloud-errors}
\end{figure}

Although SG3 yields the most accurate cloud field, its error distribution is very similar to that of OG4. Figure~\ref{fig:mswm-tc2-cloud-errors} shows the spatial error distribution for SG3 and OG4 on the uniform level-7 SCVT grid. Larger errors also appear near the equator, possibly associated with the maxima of fluid depth, temperature, and wind magnitude in the initial fields. Both schemes also show localized larger errors near pentagonal regions, which are associated with grid-imprinting effects. Such localized amplification is also seen, for example, in the fluid depth field with TRiSK, due to geometric misalignment in these regions \citep{peixoto2013gridimprinting,peixoto2016consitency}. These results indicate that, although the advection schemes are high-order accurate, they are still affected by grid imprinting, likely because the advection equation is written in conservative form using the continuity equation, which is discretized with TRiSK.

In Figure~\ref{fig:mswm-tc2-ref}, we present the errors on the variable-resolution grids.
For temperature (Fig.~\ref{fig:mswm-tc2-ref}b) and water vapor (Fig.~\ref{fig:mswm-tc2-ref}c), the SG3 scheme consistently yields smaller errors, with the improvement becoming more pronounced at higher resolutions. A similar trend, although less marked, is observed for the cloud variable (Fig.~\ref{fig:mswm-tc2-ref}d).

For temperature and water vapour, the SG2 and SG4 schemes are more accurate than the OG schemes on the coarser grids, but their errors approach those of the OG schemes at the finest resolutions.
For the cloud variable, all SG schemes exhibit similar errors, which are consistently smaller than those of the OG schemes.

In the uniform grid simulations for this test case, we observed that OG4 could overall provide smaller errors for temperature and water vapour, especially at coarser resolutions (Figure \ref{fig:mswm-tc2}). However, this behavior is not observed on the variable-resolution grids, where the advantage of OG4 diminishes.
Finally, we also observed that, overall, the errors on the variable-resolution grid follow a grid-imprinting pattern for both OG and SG schemes, particularly in regions with misaligned cells (not shown), similarly to what was reported in \citet{santos2021andes}.

\subsubsection{Barotropic instability zonal jet test case}
\label{sec:num-exp-mswm-barotropic}

We next consider the barotropically unstable zonal jet in the Southern Hemisphere from \citet{santos2021andes}. In this test, a steady zonal jet is prescribed and a small perturbation is added to trigger instability. The initial conditions are given in detail by \citet{santos2021andes}. This setup adapts the barotropic instability test of \citet{galewsky2004barotropic}, later extended to the moist shallow-water model by \citet{ferguson2019amr}.
Simulations are run for 7 days on a quasi-uniform level-7 SCVT grid with a 30-second time step, and on a topography-based refined level-7 SCVT grid with a 15-second time step.
Unlike in the previous test, we also apply fourth-order hyperdiffusion on the quasi-uniform grid to stabilize the model.

\begin{figure}[!htb]
	\centering
    \includegraphics[width=1\linewidth]{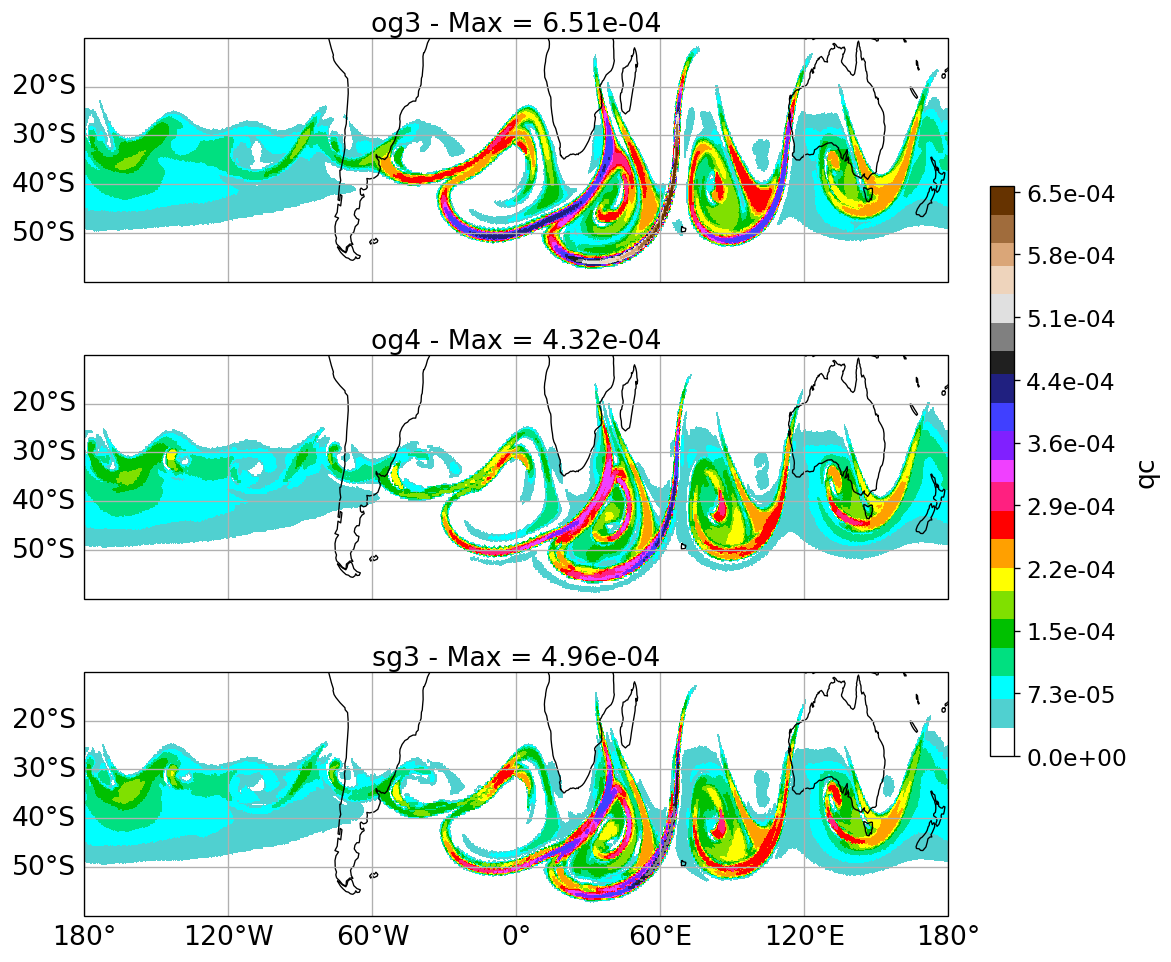}
\caption{Cloud field in the Southern Hemisphere at day 7 for the moist shallow-water model considering the barotropic instability test case, using the OG3 (top), OG4 (middle), and SG3 (bottom) advection schemes on a quasi-uniform SCVT grid at level 7.}
\label{fig:mswm-baro-unif7-cloud}
\end{figure}

\begin{figure}[!htb]
	\centering
    \includegraphics[width=1\linewidth]{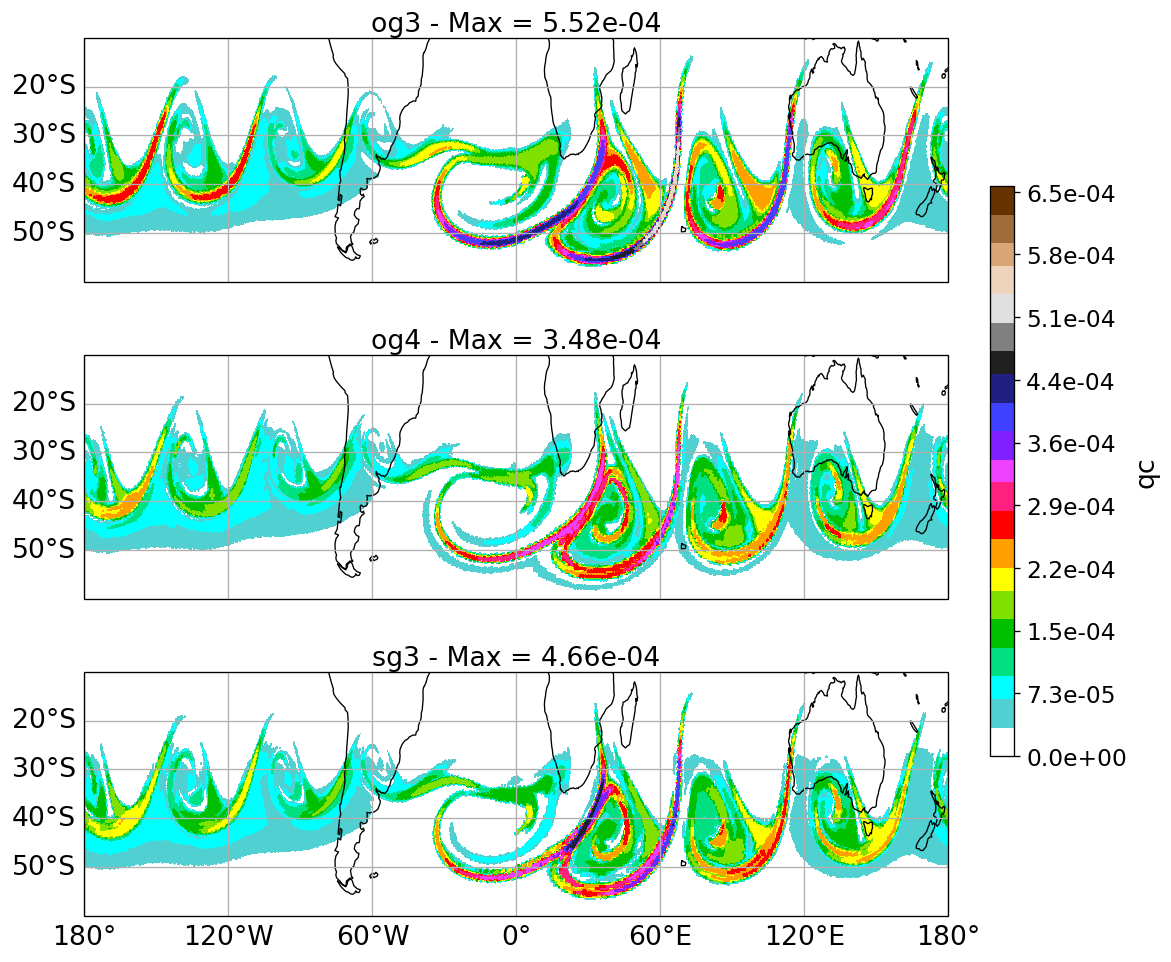}
\caption{As in Fig.~\ref{fig:mswm-baro-unif7-cloud}, but for the Andes topography–based refined SCVT grid.}
\label{fig:mswm-baro-ref7-cloud}
\end{figure}

In Figures~\ref{fig:mswm-baro-unif7-cloud} and \ref{fig:mswm-baro-ref7-cloud}, we show the cloud field after 7 days of simulation on the quasi-uniform and refined SCVT grids, respectively, for the OG3, OG4, and SG3 schemes.
All schemes capture the vortex formation and filamentary structures associated with the barotropic instability.
OG3 produces the largest cloud maxima and sharper gradients, whereas OG4 appears more diffusive, yielding smoother structures and lower peak values. SG3 exhibits intermediate behavior.
In the South American region, all schemes produce more detailed vortex structures on the Andes-refined grid than on the quasi-uniform SCVT grid.

\section{Conclusions}
\label{sec:conclusions}

In this work, we introduced a new high-order finite-volume advection scheme for tracer transport on SCVT grids. The proposed method extends the $k$-exact reconstruction approach of \citep{ollivier1997high,olliviergooch2002advdifeq} from planar unstructured meshes to the spherical surface by employing local tangent-plane projections, on which the polynomial reconstruction is performed before evaluating fluxes across cell faces. The scheme has several important characteristics: it consistently employs upwind-biased reconstructions to enhance stability, and the reconstruction function is derived by imposing local mass conservation over neighboring control volumes. We also showed that the flux limiter proposed by \citet{zalesak1979fct} can be incorporated into the scheme.

We constructed schemes of order 2, 3, and 4 using this approach,
by employing polynomial reconstructions of increasing degree and a
corresponding increase in the number of quadrature points. These schemes
were compared with the second-, third- and fourth-order schemes of
SG2011 \citep{skamarockgassmann2011}. The
comparisons were carried out on both quasi-uniform SCVT meshes and
locally refined grids generated based on topography.

The numerical experiments included standard passive advection test cases
on the sphere following the benchmark suite of \citet{nair2010class}. The
results show that the proposed second- and fourth-order schemes yield
smaller errors than their SG2011 counterparts. The third-order scheme of
SG2011, which is the only scheme in that work employing an upwind bias,
performs better than our third-order scheme. Overall, both the proposed
schemes and the third-order scheme of SG2011 show little sensitivity to
locally refined grids, whereas the second- and fourth-order schemes are
more affected by grid distortion. This behavior highlights the
stabilizing effect of the upwind bias in the third-order scheme of
SG2011. 

In general, the proposed fourth-order scheme yields the smallest errors in
the advection simulations, followed by the third-order scheme of SG2011.
However, the increased accuracy comes at a higher computational cost,
since the method employs two quadrature points per face together with a
fourth-order polynomial reconstruction. In contrast, the SG2011 third order scheme
uses a single quadrature point for the flux evaluation and a
second-order reconstruction.

We also performed simulations using the moist shallow-water model of \citet{zerroukat2015moistSWE}. In this framework, the new advection scheme was applied to the transport of water vapor, cloud water, and rainwater tracers, while the remaining dynamics were solved using the mimetic finite-volume/difference TRiSK discretization \citep{thuburn2009TRiSK, ringler2010trsk}.

One of our main findings is that grid imprinting is still observed in the tracer fields, even when high-order advection schemes are employed. This behavior is likely associated with the low-order TRiSK discretization used for the fluid depth, which appears in the flux formulation of the tracer transport equations.

A natural extension of this work would be to incorporate the proposed
advection fluxes into the momentum equations, in a manner similar to
\citet{gassmann2022momentum}. Such an extension could lead to a
higher-order shallow-water solver and provide greater consistency
between tracer advection and the dynamical core in the moist
shallow-water model. It could also help mitigate the grid imprinting
observed in the tracer fields in our simulations, despite the use of
high-order advection schemes.

\section*{Declaration of competing interest}
The authors declare that they have no known competing financial interests or personal relationships that could have appeared to influence the work reported in this paper.


\section*{Data availability}
The source code used in this work is available at \url{https://github.com/pedrospeixoto/iModel}.

\section*{Acknowledgments}
The authors acknowledge financial support from the National Council for Scientific and Technological Development (CNPq, Brazil), grant  303436/2022-0, and from the Coordenação de Aperfeiçoamento de Pessoal de Nível Superior - Brasil (CAPES), Finance Code 001. LS and PP acknowledge support from the São Paulo Research Foundation (FAPESP), grants 20/10280-4, 2021/06176-0 and 2025/02678-1.


\appendix

\section{Time Integration}
\label{sec:time-integration}
We consider the ordinary differential equation
\begin{equation}
    \frac{dX}{dt}(t) = F(X(t)), \quad X(0) = X_0,
\end{equation}
where $X$ denotes the discrete solution variables at the grid points, and $F$ is the discrete spatial operator. We advance the solution using the following three-stage Runge–Kutta scheme \citep{skamarockgassmann2011, wicker2002timeSplitting, wang2009wrfadvevaluation}:
\begin{align}
    X^{n+\frac{1}{3}} &= X^n - \frac{\Delta t}{3} 
    F(X^n), \\
    X^{n+\frac{1}{2}} &= X^n - \frac{\Delta t}{2}
    F(X^{n+\frac{1}{3}}), \\
    X^{n+1} &= X^n - \Delta t
    F(X^{n+\frac{1}{2}}).
\end{align}

\section{Flux limiter}
\label{sec:flux-limiter}

We adopt the flux-corrected transport (FCT) scheme of \citeaynum{zalesak1979fct} to limit fluxes and preserve monotonicity. As in the SG2011 method, our implementation follows the formulation of \citep{wang2009wrfadvevaluation}, extended here from Cartesian to unstructured grids.
We begin with the three-stage Runge-Kutta scheme \citep{skamarockgassmann2011}:
\begin{align}
    \label{rk3-1}
    \phi^{n+\frac{1}{3}}_i &= \phi^n_i -\frac{\Delta t}{3} 
    \sum_{e\in EC(i)} n_{e,i}F_{e}^{H}(\phi^n), \\
    \label{rk3-2}
    \phi^{n+\frac{1}{2}}_i &= \phi^n_i -\frac{\Delta t}{2}
    \sum_{e\in EC(i)} n_{e,i}F_{e}^{H}(\phi^{n+\frac{1}{3}}), \\
    \label{rk3-3}
    \phi^{n+1}_i &= \phi^n_i - \Delta t
    \sum_{e\in EC(i)} n_{e,i}F_{e}^{H}(\phi^{n+\frac{1}{2}}),
\end{align}
where $F_e^{H}$ denotes the high-order flux, e.g., the SG flux or the OG flux.
To enforce monotonicity, we modify the final stage (Equation~\eqref{rk3-3}) by first computing a low-order upwind solution:
\begin{equation}
\label{low-order-phi}
    \phi_i^{L} = \phi^n_i -{\Delta t}
    \sum_{e\in EC(i)} n_{e,i}F_{e}^{L}(\phi^n),
\end{equation}
where $F_{e}^{L}$ denotes the low-order flux. We adopt the first-order upwind flux, namely:
\begin{equation}
F_{e}^L = 
\begin{cases}
\phi_i^n, & \text{if } n_{e,i}u_{e}^n > 0, \\
\phi_j^n, & \text{otherwise},
\end{cases}
\end{equation}
with $i,j\in CE(e)$.
Next, we define the correction flux as:
\begin{equation}
     F_{e}^{C} = F_{e}^{H}(\phi^{n+\frac{1}{2}}) - F_{e}^{L}(\phi^n).
\end{equation}
We then add the flux correction to the low-order solution (Eq. \eqref{low-order-phi}), considering the outgoing and incoming fluxes separately, respectively:
\begin{align}
    \phi_i^{+} &= \phi^{L}_i -{\Delta t}
    \sum_{e\in EC(i)} \max{(n_{e,i}F_{e}^{C},0)},\\
    \phi_i^{-} &= \phi^{L}_i -{\Delta t}
    \sum_{e\in EC(i)} \min{(n_{e,i}F_{e}^{C},0)}.
\end{align}
We then define a normalization factor:
\begin{equation}
R_{e} = 
\begin{cases}
\min{
\bigg(1,
\frac{\phi^{L}_i-\phi^{min}_i}{\phi^{L}_i-\phi^{+}_i},
\frac{\phi^{L}_j-\phi^{max}_j}{\phi^{L}_j-\phi^{-}_j}
\bigg)
}, & \text{if } n_{e,i}F_{e}^{C} > 0, \\
\min{
\bigg(1,
\frac{\phi^{L}_j-\phi^{min}_j}{\phi^{L}_j-\phi^{+}_j},
\frac{\phi^{L}_i-\phi^{max}_i}{\phi^{L}_i-\phi^{-}_i}
\bigg)
}, & \text{otherwise},
\end{cases}
\end{equation}
with $i,j \in CE(e)$, $i\neq j$, $\phi^{\min}$ and $\phi^{\max}$ denote the local minimum and maximum bounds, respectively, computed over the upwind neighborhood, including the central cell, associated with the $j$-th edge of cell $i$, at time level $n$.

Finally, the last step of the Runge–Kutta scheme (Equation~\eqref{rk3-3}) is replaced by
\begin{equation}
\phi^{n+1}_i = \phi^L_i - \Delta t
\sum_{e\in EC(i)} n_{e,i}R_{e}F_{e}^{C},
\end{equation}
which can be re-expressed as
\begin{equation}
\phi^{n+1}_i = \phi^n_i - \Delta t
\sum_{e\in EC(i)} n_{e,i}\big[(1-R_e)F_{e}^{L} + R_{e}F_{e}^{H}\big],
\end{equation}
representing the limited flux as a convex combination of the low-order and high-order fluxes.

\section{Wind velocity reconstruction at quadrature points}
\label{sec:wind-recon}
\begin{figure}[!htb]
    \centering
    \includegraphics[width=0.8\textwidth]{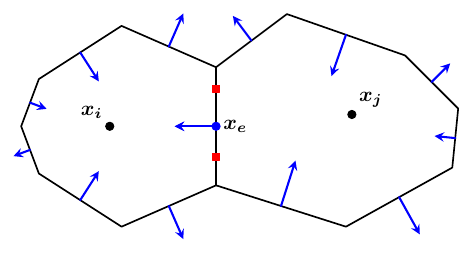}
    \caption{Normal velocities at edge points (blue arrows) and quadrature points on edges (red squares) used for the local reconstruction centered at $\boldsymbol{x}_e$ (blue circle).}
\label{fig:normal_wind}
\end{figure}

The momentum equation \eqref{mswm-ode-mom} provides the normal velocity at edge points $\boldsymbol{x}_e$, which is sufficient for the SG scheme. The OG scheme, however, requires velocities at quadrature points.
Given the normal velocities $u_\ell$ at the edge points $\boldsymbol{x}_{\ell}$, shown by the blue arrows in Fig.~\ref{fig:normal_wind}, we reconstruct the normal velocity at the quadrature points on each edge $e$, shown by the red squares.

As in the OG and SG schemes, we perform the reconstruction on a tangent plane centered at $\boldsymbol{x}_e$, onto which we project the edge points $\boldsymbol{x}_{\ell}$ for $\ell \in EC(i) \cup EC(j)$, with $i,j \in CE(e)$. The reconstruction uses the normal velocities on the edges of cells $i$ and $j$.
Following \citeaynum{peixoto2014vecrecon}, we seek a linear approximation of the velocity field of the form
\begin{equation}
\mathbf{u}(x,y) = \mathbf{a}_0 + \mathbf{a}_1 x + \mathbf{a}_2 y,
\end{equation}
where $(x,y)$ are the local tangent plane coordinates, $\mathbf{a}_0 = (a_0^x,a_0^y)$, $\mathbf{a}_1 = (a_1^x,a_1^y)$, and $\mathbf{a}_2 = (a_2^x,a_2^y)$.

For $\ell \in EC(i) \cup EC(j)$, we impose
\begin{equation}
\mathbf{u}(x_\ell,y_\ell)\cdot \boldsymbol{\hat{n}}_\ell
= \mathbf{a}_0 \cdot \boldsymbol{\hat{n}}_\ell 
+ (\mathbf{a}_1\cdot \boldsymbol{\hat{n}}_\ell)\, x_\ell 
+ (\mathbf{a}_2\cdot \boldsymbol{\hat{n}}_\ell)\, y_\ell
= u_\ell,
\end{equation}
where $\boldsymbol{\hat{n}}_\ell$ is the normal vector projected onto the tangent plane. This defines an overdetermined system, solved in the least-squares sense for $(a_0^x,a_0^y,a_1^x,a_1^y,a_2^x,a_2^y)$, yielding a second-order velocity reconstruction \citep{peixoto2014vecrecon}. We also include distance-based weights, since nearby edge points provide more accurate information.




\bibliographystyle{elsarticle-num-names}
\bibliography{cas-refs}

\end{document}